# A Comparative Study of Numerical Methods for Approximating the Solutions of a Macroscopic Automated-Vehicle Traffic Flow Model

**George Titakis[*], Iasson Karafyllis[**], Dionysis Theodosis[*], Ioannis Papamichail[*] and Markos Papageorgiou[*,***]**

[*] Dynamic Systems and Simulation Laboratory,
Technical University of Crete, Chania, 73100, Greece,

emails: gtitakis@dssl.tuc.gr, dtheodosis@dssl.tuc.gr,
ipapa@dssl.tuc.gr, markos@dssl.tuc.gr

[**] Dept. of Mathematics, National Technical University of Athens,
Zografou Campus, 15780, Athens, Greece,
emails: iasonkar@central.ntua.gr , iasonkaraf@gmail.com

[***] Faculty of Maritime and Transportation,
Ningbo University, Ningbo, China.

**Abstract**

In this paper, a particle method is used to approximate the solutions of a "fluid-like" macroscopic traffic flow model for automated vehicles. It is shown that this method preserves certain differential inequalities that hold for the macroscopic traffic model: mass is preserved, the mechanical energy is decaying and an energy functional is also decaying. To demonstrate the advantages of the particle method under consideration, a comparison with other numerical methods for viscous compressible fluid models is provided. Since the solutions of the macroscopic traffic model can be approximated by the solutions of a reduced model consisting of a single nonlinear heat-type partial differential equation, the numerical solutions produced by the particle method are also compared with the numerical solutions of the reduced model. Finally, a traffic simulation scenario and a comparison with the Aw-Rascle-Zhang (ARZ) model are provided, illustrating the advantages of the use of automated vehicles.

## 1. Introduction

Traffic flow theory is based on macroscopic models to describe the aggregate behaviour of the vehicles in a traffic stream. In these models, traffic flow may be viewed as a fluid with specific characteristics expressed in aggregate variables like vehicle density, flow, and mean speed of vehicles, while its evolution is described via Partial Differential Equations (PDEs). Depending on the number of the PDEs, macroscopic models are categorized as first-order and second-order models. Among the first-order models for human drivers, the most well-known was independently proposed by Lighthill and Whitham [25] and Richards [30] (LWR model). As a first-order model, the LWR model consists of a single hyperbolic PDE that describes the mass (vehicle) conservation law. However, the inadequacy of first-order models to reproduce some real dynamic phenomena led to the development of second-order models, such as the well-known ARZ (Aw-Rascle-Zhang) model (see [3], [35]), which is supplemented with an additional PDE that expresses the dynamics of mean speed.

In contrast to traffic flow models for human drivers, the macroscopic traffic model derived in [17] via the design of bidirectional cruise controllers for automated vehicles (AVs) is isotropic, as in other fluids, with the principle that vehicles react to both upstream and downstream vehicles. Isotropy, as well as many other characteristics of the compressible fluid flow that are present in the derived macroscopic traffic model for AVs, follow from the fact that nudging and lane-free vehicle movement are allowed by the cruise controllers. Moreover, a link between the macroscopic model under consideration and the underlying microscopic cruise controllers is provided in [17], which allows to determine the physical properties of the "traffic-fluid" by changing the functions or the parameter values included in the cruise controllers.

In the present work, we deal with the approximation of the solutions of a macroscopic traffic model for AVs which follows from a certain selection of the potential function involved in the cruise controller (see [34]). This selection, in addition to providing more "fluid-like" characteristics, enables a model reduction. The resulting reduced model consists of a single nonlinear heat-type partial differential equation and approximates the solutions of the macroscopic AV model. This fact is analogous to the human-driven vehicles case, where the ARZ model reduces to the LWR model [24]. For the derivation of the macroscopic AV model, a particle method was proposed in [17], in which the vehicles are considered as "self-driven" particles of a fluid. This method is a modified version of the particle method presented in [16] for compressible fluids, and we use it for approximating the solutions of the macroscopic AV model. It is shown that all differential inequalities that hold for the macroscopic AV model (described in Section 2) are preserved by the particle method (see Proposition 1 and Proposition 3 below). It is noted that, to the best of our knowledge, a numerical study of the particle method and comparison with other numerical methods has never been performed in the literature for viscous compressible fluid models. Therefore, a main goal of the present work is the numerical study of the particle scheme and the derivation of important conclusions about its effectiveness.

Although many numerical methods have been used for compressible fluid models with constant viscosity (see for example [1], [2], [8], [11], [13], [21] and references therein), there are very few articles that cover the case of density-dependent viscosity (see [4], [5], [9], [10], [16], [29]) and most of them utilize an approximation scheme for proving existence of solutions. Since the macroscopic AV model that we study is similar to compressible fluids with density-dependent viscosity (see Section 2), we compare the particle method proposed in [16], [17] with the so-called "method of lines" used in [10], [29] for proving existence of solutions for compressible fluid models with density-dependent viscosity, as well as with the numerical method used in [9] for approximating the solutions of the viscous Saint-Venant equations.

A common characteristic between the method of lines ([10], [29]) and the particle method proposed in [16], [17] is that in both methods the solutions of the macroscopic AV model are approximated by solving a system of Ordinary Differential Equations (ODEs). However, with the method of lines, the solutions are approximated at fixed grid-points in Lagrangian coordinates, and, as a result, the actual coordinates of the vehicles remain unknown. Therefore, a variable transformation is required to obtain the solutions in Eulerian coordinates, which, in addition to having an impact on the computational cost, (see Section 4), may create an additional error in the approximation of the solutions. On the other hand, the approximation of the solutions with the particle method proposed in [16] is given directly in Eulerian coordinates, and the computational cost, in terms of CPU time, is much less compared to the method of lines and the numerical method used in [9]. To demonstrate the advantages of the particle method as well as the properties of the macroscopic AV model, two examples are presented. The first example focuses on evaluating the performance of the proposed particle method, while the second example includes also a realistic traffic simulation scenario for comparing the macroscopic AV model with the ARZ model. In both examples, the particle method is compared with the method of lines ([10], [29]) and the numerical method used for the viscous Saint-Venant equations [9] showing that, although the difference between the numerical solutions is small, higher-resolution approximations are provided by the proposed particle method. Comparisons with the numerical solution of the reduced model (see [34]) are also provided for both examples, indicating its high-resolution approximations.

The structure of the paper is as follows. Section 2 presents the macroscopic traffic model for AVs as well as the particle method used for approximating its solutions. Section 3 is devoted to the presentation of some numerical methods used in the literature for approximating the solutions of viscous compressible fluid models. Section 4 presents the two aforementioned examples. Finally, Section 5 gives the concluding remarks of the present work. All proofs are provided in the Appendix.

**Notation.** Throughout this paper, we adopt the following notation.

* By $|x|$ we denote both the Euclidean norm of a vector $x \in \Re^n$ and the absolute value of a scalar $x \in \Re$.

* $\Re_+ := [0,+\infty)$ denotes the set of non-negative real numbers.

* By $x'$ we denote the transpose of a vector $x \in \Re^n$. By $|x|_\infty = \max\{|x_i|, i = 1,...,n\}$ we denote the infinity norm of a vector $x = (x_1, x_2,...,x_n)' \in \Re^n$.

* Let $A \subseteq \Re^n$ be an open set. By $C^0(A;\Omega)$, we denote the class of continuous functions on $A \subseteq \Re^n$, which take values in $\Omega \subseteq \Re^m$. By $C^k(A;\Omega)$, where $k \geq 1$ is an integer, we denote the class of functions on $A \subseteq \Re^n$ with continuous derivatives of order $k$, which take values in $\Omega \subseteq \Re^m$. When $\Omega = \Re$ then we write $C^0(A)$ or $C^k(A)$.

* Let $I \subseteq \Re$ be a given interval. For $p \in [1,\infty)$, $L^p(I)$ denotes the set of equivalence classes of Lebesgue measurable functions $f : I \to \Re$ with $\|f\|_p := \left(\int_I |f(x)|^p dx\right)^{1/p} < +\infty$. $L^\infty(I)$ denotes the set of equivalence classes of measurable functions $f : I \to \Re$ for which $\|f\|_\infty = \text{ess}\sup_{x \in I}(|f(x)|) < +\infty$.

* Let $u : \left(\bigcup_{t \geq 0} \{t\} \times I(t)\right) \to \Re$ where $I(t) \subseteq \Re$ is an open interval for each $t \geq 0$, be any function differentiable with respect to its arguments. We use the notation $u_t(t,x) = \frac{\partial u}{\partial t}(t,x)$ and $u_x(t,x) = \frac{\partial u}{\partial x}(t,x)$ for the partial derivatives of $u$ with respect to $t$ and $x$, respectively. We use the notation $u[t]$ to denote the profile at certain $t \geq 0$

, $(u[t])[x] := u(t,x)$, for all $x \in I(t)$. We also use the notation $u(t,x^-)$ and $u(t,x^+)$ for the limits $u(t,x^-) = \lim_{s \to x^-}(u(t,s))$ and $u(t,x^+) = \lim_{s \to x^+}(u(t,s))$.

## 2. A Particle Method for a Macroscopic Automated-Vehicle Traffic Flow Model

*2.1 The Macroscopic Traffic Flow Model for Automated Vehicles*

The macroscopic traffic flow model for automated vehicles (2.1)-(2.4) which arises from the study of microscopic vehicle movement control laws (cruise controllers) for automated vehicles on lane-free roads [17], [34] is described as follows for all $\tau > 0$, $\xi \in (\tilde{l}(\tau), \tilde{L}(\tau))$:

$$\tilde{\rho}_\tau + (\tilde{\rho}\tilde{v})_\xi = 0 \tag{2.1}$$

$$\tilde{\rho}\tilde{\beta}'(\tilde{v})\tilde{v}_\tau + \tilde{\rho}\tilde{\beta}'(\tilde{v})\tilde{v}\tilde{v}_\xi + \tilde{\sigma}\tilde{\kappa}(\tilde{\rho})\tilde{\rho}_\xi = (\tilde{\rho}\tilde{\kappa}(\tilde{\rho})\tilde{v}_\xi)_\xi - \tilde{\sigma}\tilde{\rho}\tilde{\beta}(\tilde{v}) \tag{2.2}$$

$$\frac{d}{d\tau}\tilde{l}(\tau) = \tilde{v}\left(\tau, \tilde{l}^+(\tau)\right) \tag{2.3}$$

$$\frac{d}{d\tau}\tilde{L}(\tau) = \tilde{v}\left(\tau, \tilde{L}^-(\tau)\right) \tag{2.4}$$

The states $\tilde{l}(\tau), \tilde{L}(\tau) \in \Re$ yield the interval $(\tilde{l}(\tau), \tilde{L}(\tau))$ that contains the vehicles for times $\tau > 0$; $\tilde{\rho}(\tau, \xi) \in (0, \rho_{max})$, $\tilde{v}(\tau, \xi) \in (0, v_{max})$ are the traffic density and mean speed, respectively, at position $\xi \in (\tilde{l}(\tau), \tilde{L}(\tau))$ on a highway, with $\rho_{max} > 0$, $v_{max} > 0$ being the maximum density and speed limits, respectively; $v^* \in (0, v_{max})$ is the speed set-point and $\tilde{\sigma} > 0$ is a constant. Model (2.1)-(2.4) follows from a specific selection of the functions included in the cruise controller (see [34]), inspired by isentropic compressible fluid flow models (see [26], [33]). The function $\tilde{\beta}(\tilde{v})$ is defined as follows for all $\tilde{v} \in (0, v_{max})$:

$$\tilde{\beta}(\tilde{v}) = \int_{v^*}^{\tilde{v}} \tilde{q}(s)ds \tag{2.5}$$

where

$$\tilde{q}(\tilde{v}) = \frac{v_{max}^3(\tilde{v}+v^*) - 2v_{max}^2 \tilde{v}v^*}{2(v_{max}-\tilde{v})^2 \tilde{v}^2} \text{ for all } \tilde{v} \in (0, v_{max}). \tag{2.6}$$

In addition, the function $\tilde{\kappa}:(0, \rho_{max}) \to \Re_+$, is $C^1((0, \rho_{max}))$ and satisfies the following properties

$$\lim_{\tilde{\rho} \to \rho_{max}^-} \tilde{\kappa}(\tilde{\rho}) = +\infty \tag{2.7}$$

$$\tilde{\kappa}(\tilde{\rho}) = 0 \text{ for all } \tilde{\rho} \in (0, \bar{\rho}] \text{ and } \tilde{\kappa}(\tilde{\rho}) > 0 \text{ for all } \tilde{\rho} \in (\bar{\rho}, \rho_{max}) \tag{2.8}$$

where the constant $\bar{\rho} \in (0, \rho_{max})$ is called "the interaction density" since for $\tilde{\rho} \in (0, \bar{\rho}]$ there is no interaction among vehicles (see (2.8)). It should be noted that outside the interval $(\tilde{l}(\tau), \tilde{L}(\tau))$ the highway is "empty" and therefore:

$$\tilde{\rho}(\tau, \xi) = 0 \text{ for } \xi \leq \tilde{l}(\tau) \text{ and } \xi \geq \tilde{L}(\tau), \tau > 0. \quad (2.9)$$

Furthermore, it is possible (see [34]) to approximate the solutions of (2.1)-(2.4) by the solutions of a reduced model that consists of the following single nonlinear heat-type Partial Differential Equation (PDE):

$$\tilde{\rho}_\tau + (\tilde{\rho}\tilde{\beta}^{-1}(-\tilde{\rho}^{-1}\tilde{\kappa}(\tilde{\rho})\tilde{\rho}_\xi))_\xi = 0. \quad (2.10)$$

Moreover, notice that the speed PDE (2.2) is strongly reminiscent of the momentum equation for 1-D polytropic compressible fluid flow since it contains a viscosity term with density-dependent viscosity $(\tilde{\rho}\tilde{\kappa}(\tilde{\rho})\tilde{v}_\xi)_\xi$, a pressure term $\tilde{\sigma}\tilde{\kappa}(\tilde{\rho})\tilde{\rho}_\xi$ and a friction term $-\tilde{\sigma}\tilde{\rho}\tilde{\beta}(\tilde{v})$ (by analogy with the relaxation term involved in the macroscopic traffic flow models for human drivers, [3], [35]). The pressure term $\tilde{\sigma}\tilde{\kappa}(\tilde{\rho})\tilde{\rho}_\xi$ describes the tendency of vehicles to accelerate or decelerate based on the (local) density, the term $(\tilde{\rho}\tilde{\kappa}(\tilde{\rho})\tilde{v}_\xi)_\xi$ is a viscosity term with $\tilde{\rho}\tilde{\kappa}(\tilde{\rho})$ playing the role of viscosity and making the "traffic fluid" act as a Newtonian fluid, while the friction term $-\tilde{\sigma}\tilde{\rho}\tilde{\beta}(\tilde{v})$ expresses the tendency of vehicles to adjust their speed to the given speed set-point.

**Remarks: (i)** In contrast to conventional traffic, the traffic flow for model (2.1)-(2.4) is isotropic, as in fluid flow, since the vehicles are automated and react to both upstream and downstream vehicles.
**(ii)** The selection of $\tilde{\kappa}(\tilde{\rho})$ to satisfy condition (2.7) has the effect of making the fluid behave like a solid when the density tends to the maximum density $\rho_{max}$.

Assume that $r > 0$ is a given constant length. Using the variable transformation $\xi = rx + v^*\tau$, $x \in (l(t), L(t))$ where $l(t) = \dfrac{\tilde{l}(\tau) - v^*\tau}{r}$, $L(t) = \dfrac{\tilde{L}(\tau) - v^*\tau}{r}$ and the dimensionless quantities $\tau = \dfrac{r}{v^*}t$, $b = \dfrac{v_{max} - v^*}{v^*}$, $R = \dfrac{\rho_{max}}{\bar{\rho}} > 1$, $w = \dfrac{\tilde{v} - v^*}{v^*}$, $\rho = \dfrac{\tilde{\rho}}{\bar{\rho}}$, $\sigma = \dfrac{r}{v^*}\tilde{\sigma}$, we obtain the following dimensionless model for all $t > 0$, $x \in (l(t), L(t))$:

$$\rho_t + (\rho w)_x = 0 \quad (2.11)$$

$$\rho\beta'(w)w_t + \rho\beta'(w)ww_x + \sigma\kappa(\rho)\rho_x = (\rho\kappa(\rho)w_x)_x - \sigma\rho\beta(w) \quad (2.12)$$

$$\frac{d}{dt}l(t) = w(t, l^+(t)) \quad (2.13)$$

$$\frac{d}{dt}L(t) = w(t, L^-(t)) \quad (2.14)$$

with state-constraints $\rho(t,x) \in (0, R)$, $w(t,x) \in (-1, b)$ for $t > 0$, $x \in (l(t), L(t))$, where

$R > 1$, $b > 0$, $\sigma > 0$ are constants. The function $\beta(w)$ is defined by

$$\beta(w) = \int_0^w q(s)ds = \frac{b+1}{2}\left[\frac{w(b+1)}{(w+1)(b-w)} + \ln\left(\frac{b(w+1)}{b-w}\right)\right], \text{ for all } w \in (-1, b) \quad (2.15)$$

where

$$q(w) = (1+b)^2 \frac{2b + (b-1)w}{2(b-w)^2(1+w)^2}, \text{ for all } w \in (-1, b) \quad (2.16)$$

while the function $\kappa : (0, R) \to \Re_+$ is a function of class $C^1((0, R))$ that satisfies the following properties:

$$\lim_{\rho \to R^-} \kappa(\rho) = +\infty \quad (2.17)$$

$$\kappa(\rho) = 0 \text{ for all } \rho \in (0, 1] \text{ and } \kappa(\rho) > 0 \text{ for all } \rho \in (1, R). \quad (2.18)$$

Notice that outside the interval $(l(t), L(t))$ it holds that:

$$\rho(t, x) = 0 \text{ for } x \leq l(t) \text{ and } x \geq L(t), t > 0. \quad (2.19)$$

The solutions of (2.11)-(2.14) can be approximated by the solutions of the following dimensionless reduced model:

$$\rho_t + (\rho \beta^{-1}(-\rho^{-1}\kappa(\rho)\rho_x))_x = 0. \quad (2.20)$$

The relations between the functions $q, \beta, \kappa$ of the dimensionless model (2.11)-(2.14) and the functions $\tilde{q}, \tilde{\beta}, \tilde{\kappa}$ of the original model (2.1)-(2.4) are given in Table 1.

| Dimensionless model (2.11)-(2.14) | Original model (2.1)-(2.4) |
|---|---|
| $q(w)$ | $= \tilde{q}(v^*(w+1))$ |
| $\beta(w)$ | $= \dfrac{1}{v^*}\tilde{\beta}(v^*(w+1))$ |
| $\kappa(\rho)$ | $= \dfrac{1}{rv^*}\tilde{\kappa}(\rho\bar{\rho})$ |

**Table 1**: Relations between the functions $q, \beta, \kappa$ of the dimensionless model (2.11)-(2.14) and the functions $\tilde{q}, \tilde{\beta}, \tilde{\kappa}$ of the original model (2.1)-(2.4).

For any classical solution of model (2.11)-(2.14) we can define the following functionals:

$$m(t) = \int_{l(t)}^{L(t)} \rho(t, x)dx \quad (2.21)$$

$$E(t) = \int_{l(t)}^{L(t)} \rho(t, x) H(w(t, x))dx + \int_{l(t)}^{L(t)} Q(\rho(t, x))dx \quad (2.22)$$

$$W(t) = \frac{1}{2} \int_{l(t)}^{L(t)} \frac{1}{\rho(t,x)} \left( \rho(t,x)\beta(w(t,x)) + \sigma^{-1} \left( P(\rho(t,x)) \right)_x \right)^2 dx \quad (2.23)$$

where

$$H(w) = \int_0^w sq(s)ds = \frac{w^2(b+1)^2}{2(w+1)(b-w)}, \text{ for all } w \in (-1,b) \quad (2.24)$$

$$Q(\rho) = \rho \int_1^\rho \tau^{-2} P(\tau) d\tau, \text{ for all } \rho \in (0,R) \quad (2.25)$$

$$P(\rho) = \sigma \int_1^\rho \kappa(\tau) d\tau, \text{ for all } \rho \in (0,R). \quad (2.26)$$

Notice that definitions (2.25), (2.26) in conjunction with (2.18) imply that $Q(\rho) \geq 0$, $P(\rho) \geq 0$ for all $\rho \in (0,R)$. The functional (2.21) expresses the total mass and the functional (2.22) is inspired by the mechanical energy of model (2.11)-(2.14) where $\int_{l(t)}^{L(t)} \rho(t,x) H(w(t,x)) dx$ denotes the kinetic energy and $\int_{l(t)}^{L(t)} Q(\rho(t,x)) dx$ the potential energy. Moreover, the functional (2.23) is an energy functional that has been constructed using the following transformation

$$\varphi(t,x) = \rho(t,x)\beta(w(t,x)) + \sigma^{-1}\left(P(\rho(t,x))\right)_x. \quad (2.27)$$

The transformation (2.27) has been used in the literature of isentropic, compressible fluid flow (see [14], [15], [19], [20], [22], [31], [32], [33]) and has the effect of making the viscosity disappear from equation (2.12), which then becomes

$$\varphi_t + (w\varphi)_x + (P(\rho))_x = -\sigma\rho\beta(w). \quad (2.28)$$

A direct consequence of (2.21)-(2.26), (2.19) and equations (2.11)-(2.14) are the following estimates

$$\dot{m}(t) = 0 \quad (2.29)$$

$$\dot{E}(t) = -\int_{l(t)}^{L(t)} \rho(t,x)\kappa(\rho(t,x)) w_x^2(t,x) dx - \sigma \int_{l(t)}^{L(t)} \rho(t,x) w(t,x) \beta(w(t,x)) dx \leq 0 \quad (2.30)$$

$$\dot{W}(t) = -\sigma \int_{l(t)}^{L(t)} \frac{1}{\rho(t,x)} \left( \rho(t,x)\beta(w(t,x)) + \sigma^{-1}\left(P(\rho(t,x))\right)_x \right)^2 dx \leq 0 \quad (2.31)$$

which hold for all $t > 0$. Definition (2.29) implies that the total mass, defined by (2.21), remains constant while (2.30), (2.31) show that the functionals defined by (2.22), (2.23) are non-increasing.

## 2.2 The Particle Method

Several particle methods in the literature have been used for the numerical solution of partial differential equations (see for example [6], [7]). Here, we consider a particle method that was initially proposed in [17] to derive model (2.11)-(2.14) by designing cruise controllers on a microscopic scale. A similar particle method (without a friction term) was also recently presented [16], that provided approximations that converge to a weak solution of the Navier-Stokes equations for the 1-D flow of a viscous compressible fluid. The particle method follows by considering $n$ vehicles/particles with coordinates $L^{(n)}(t) = x_1(t) > x_2(t) > .... > x_n(t) = l^{(n)}(t)$ which move according to the following equations (recall that $b > 0$, $\sigma > 0$, $R > 1$):

$$s_i = x_{i-1} - x_i, \ i = 2,...,n \tag{2.32}$$

$$\dot{x}_i = w_i, \ i = 1,...,n \tag{2.33}$$

$$\beta'(w_1)\dot{w}_1 = -\sigma\beta(w_1) - na\Phi'(nas_2) + n^2 aK(nas_2)(w_2 - w_1) \tag{2.34}$$

$$\beta'(w_i)\dot{w}_i = -\sigma\beta(w_i) + na\Phi'(nas_i) - na\Phi'(nas_{i+1}) + n^2 aK(nas_i)(w_{i-1} - w_i)$$
$$+ n^2 aK(nas_{i+1})(w_{i+1} - w_i), \ i = 2,...,n-1 \tag{2.35}$$

$$\beta'(w_n)\dot{w}_n = -\sigma\beta(w_n) + na\Phi'(nas_n) + n^2 aK(nas_n)(w_{n-1} - w_n) \tag{2.36}$$

with state-space

$$\Omega_n = \left\{ (x_1,...,x_n, w_1,...,w_n) \in \Re^n \times (-1, b)^n : na(x_{i-1} - x_i) > 1/R, \ i = 2,...,n \right\} \tag{2.37}$$

where $a$ is a positive constant and $\beta(w)$ is defined by (2.15). The functions $\Phi : \left(\frac{1}{R}, +\infty\right) \to \Re_+$, $K : \left(\frac{1}{R}, +\infty\right) \to \Re_+$, $\Phi \in C^2$, $K \in C^2$ satisfy the following equations for all $s \in \left(\frac{1}{R}, +\infty\right)$:

$$K(s) = \frac{a}{s^2} \kappa\left(\frac{1}{s}\right) \tag{2.38}$$

$$\Phi(s) = sQ\left(\frac{1}{s}\right). \tag{2.39}$$

It follows from (2.18), (2.25), (2.26), (2.38), and (2.39) that

$$K(s), \Phi(s) > 0 \text{ for all } s \in \left(\frac{1}{R}, 1\right) \text{ and } K(s) = \Phi(s) = 0 \text{ for all } s \geq 1. \tag{2.40}$$

Moreover, by virtue of (2.39), (2.25) and (2.26) we get that:

$$\lim_{d \to (1/R)^+} (\Phi(d)) = +\infty \tag{2.41}$$

and

$$\Phi'(s) = -\sigma \int_1^{1/s} \kappa(\tau) d\tau \text{ for all } s \in \left(\frac{1}{R}, +\infty\right). \tag{2.42}$$

Notice that equations (2.38), (2.42) provide a link between the function $\kappa(\rho)$ used in the macroscopic model (2.11)-(2.14) and the functions $K(s)$, $\Phi'(s)$ used in the microscopic model (2.32)-(2.36). This fact implies that by changing the functions $K(s), \Phi(s)$ we can determine the physical properties of the "traffic fluid" through the function $\kappa(\rho)$.

We next define the following functions:

$$E_n(t) = \frac{1}{na} \sum_{i=1}^{n} H(w_i(t)) + \frac{1}{na} \sum_{i=2}^{n} \Phi(nas_i(t)) \tag{2.43}$$

$$W_n(t) = \frac{1}{2na} \sum_{i=2}^{n-1} \varphi_i^2(t) \tag{2.44}$$

where

$$\varphi_i(t) = \beta(w_i(t)) + \sigma^{-1} na \left(\Phi'(nas_{i+1}(t)) - \Phi'(nas_i(t))\right). \tag{2.45}$$

For the parameterized family of finite-dimensional systems (2.32)-(2.36) defined on $\Omega_n$ (recall (2.37)) with initial conditions $(x_1(0),...,x_n(0), w_1(0),...,w_n(0))$ in $\Omega_n$, the following propositions hold.

**Proposition 1:** *For every solution of (2.32)-(2.36) the following equation holds*

$$\dot{E}_n(t) = -n \sum_{i=2}^{n} K(nas_i(t))(w_{i-1}(t) - w_i(t))^2 - \frac{\sigma}{na} \sum_{i=1}^{n} w_i(t) \beta(w_i(t)) \leq 0 \tag{2.46}$$

*as long as the solution of (2.32)-(2.36) is defined.*

**Remark:** Inequality (2.46) is the discretized version of inequality (2.30) and shows that the discretized mechanical energy $E_n$ is not increasing.

**Proposition 2:** *Every solution of (2.32)-(2.36) is defined for all $t \geq 0$.*

**Proposition 3:** *For every solution of (2.32)-(2.36) the following equation holds for all $t \geq 0$:*

$$\dot{W}_n(t) = -\frac{\sigma}{na} \sum_{i=2}^{n-1} \varphi_i^2(t) \leq 0. \tag{2.47}$$

**Remark:** Inequality (2.47) is the discretized version of inequality (2.31) and shows that the discretized energy functional $W_n$ is not increasing.

Let a solution $(x_i(t), w_i(t))$, $i = 1,...,n$ of (2.32)-(2.36), defined for $t \geq 0$, be given. We define for all $t \geq 0$:

$$\rho_i(t) = \frac{1}{na(x_{i-1}(t) - x_i(t))}, \ i = 2,...,n \tag{2.48}$$

$$\rho_1(t) = \rho_2(t) = \frac{1}{na(x_1(t) - x_2(t))} \tag{2.49}$$

$$\rho^{(n)}(t,x) = \rho_i(t) + (\rho_{i-1}(t) - \rho_i(t))\left(\frac{x - x_i(t)}{x_{i-1}(t) - x_i(t)}\right), \tag{2.50}$$
$$x \in [x_i(t), x_{i-1}(t)), \ i = 2,...,n-1$$

$$w^{(n)}(t,x) = w_i(t) + (w_{i-1}(t) - w_i(t))\left(\frac{x - x_i(t)}{x_{i-1}(t) - x_i(t)}\right), \tag{2.51}$$
$$x \in [x_i(t), x_{i-1}(t)), \ i = 2,...,n-1$$

$$\rho^{(n)}(t,x) = \rho_n(t) + (\rho_{n-1}(t) - \rho_n(t))\left(\frac{x - x_n(t)}{x_{n-1}(t) - x_n(t)}\right), \ x \in (x_n(t), x_{n-1}(t)) \tag{2.52}$$

$$w^{(n)}(t,x) = w_n(t) + (w_{n-1}(t) - w_n(t))\left(\frac{x - x_n(t)}{x_{n-1}(t) - x_n(t)}\right), \ x \in (x_n(t), x_{n-1}(t)) \tag{2.53}$$

$$\begin{aligned} l^{(n)}(t) &= x_n(t) \\ L^{(n)}(t) &= x_1(t) \end{aligned} \tag{2.54}$$

The functions $\rho^{(n)} : \bigcup_{t \geq 0} \{t\} \times (x_n(t), x_1(t)) \to (0, R)$, $w^{(n)} : \bigcup_{t \geq 0} \{t\} \times (x_n(t), x_1(t)) \to (-1, b)$ are $C^1$ on the set $\bigcup_{t \geq 0} \{(t, x) : x \in (x_n(t), x_1(t)), \ x \neq x_i(t), \ i = 2,...,n-1\}$ and can be used to approximate the solutions of (2.11)-(2.14). In [16], it was proved that similar approximations converge to a weak solution of the Navier-Stokes equations for the 1-D flow of a viscous compressible fluid.

## 3. Numerical Methods for Viscous Compressible Fluids

To evaluate the performance of the proposed particle method (2.32)-(2.37), it is of great importance to compare it with other numerical methods in the literature that are used for the numerical approximation of viscous compressible fluid models with density-dependent viscosity. Equation (2.12) can be written in the following form for all $t > 0$ and $x \in (l(t), L(t))$:

$$(\rho\beta(w))_t + (\rho w \beta(w) + P(\rho))_x = (\rho\kappa(\rho)w_x)_x - \sigma\rho\beta(w) \tag{3.1}$$

where $\beta(w)$, $P(\rho)$ are defined by (2.15), (2.26), respectively. Equation (3.1) is similar to the momentum equation of the 1-D compressible Navier-Stokes equations with density-dependent viscosity (see [10], [16]). Therefore, it is interesting to compare the particle method under consideration with the so-called "method of lines" used in [10] for proving existence of global weak solutions to the 1-D compressible Navier-Stokes equations, as well as with the numerical method used in [9] for the numerical approximation of the viscous Saint-Venant equations. However, there are crucial differences between the momentum equations used in [9], [10] and equation (3.1) related to the functions $\beta(w)$, $P(\rho)$, $\mu(\rho)$ and the friction term (see Table 2). Therefore, appropriate modifications should be made to compare the proposed particle method with the aforementioned numerical methods. Note that apart from the differences described in Table 2, the momentum equation used in [9] includes also an additional term $g\rho z_x$, where $z$ is the topography of the fluid and $g$ is the acceleration of gravity (used also in the pressure function (see Table 2)).

|  | **Equation (3.1)** | **Momentum equation used in [9]** | **Momentum equation used in [10]** |
|---|---|---|---|
|  | $\beta(w) = \dfrac{b+1}{2}\left[\dfrac{w(b+1)}{(w+1)(b-w)} + \ln\left(\dfrac{b(w+1)}{b-w}\right)\right]$ | $\beta(w) = w$ | $\beta(w) = w$ |
| **Pressure** | $P(\rho) = \sigma\int_1^\rho \kappa(\tau)d\tau$ <br> $\sigma > 0$ | $P(\rho) = \sigma\rho^2$ <br> $\sigma = \dfrac{1}{2}g,\ g > 0$ | $P(\rho) = \sigma\rho^\gamma,\ \gamma > 1$ <br> $\sigma = 1$ |
| **Viscosity** | $\mu(\rho) = \rho\kappa(\rho)$ | $\mu(\rho) = \mu\rho$ <br> $\mu > 0$ | $\mu(\rho) = \mu(\theta\rho^\theta + \varepsilon\eta\rho^\eta)$ <br> $\theta, \varepsilon > 0,\ \eta \in (0,1)$ <br> $\mu = 1$ |
| **Friction Term** | $-\sigma\rho\beta(w)$ | $-\dfrac{\kappa}{1+\dfrac{\kappa\rho}{3\mu}}w$ <br> $\kappa, \mu > 0$ | No friction term |

**Table 2**: Differences between equation (3.1) and the momentum equations used in [9], [10].

To obtain a numerical approximation of the solutions of model (2.11)-(2.14), we modify the numerical method used in [9] by omitting the additional term $g\rho z_x$ from the momentum equation and substituting the functions $\beta(w)$, $P(\rho)$, $\mu(\rho)$ as well as the friction term with the corresponding functions of equation (3.1) (see Table 2). Thus, the modified version of the numerical method discussed in [9] reads as follows:

$$\rho_i^0 = \frac{1}{\delta x}\int_{x_{i-1/2}}^{x_{i+1/2}} \rho(0,x)dx,\ i \in \mathbb{Z} \qquad (3.2)$$

$$w_{i+1/2}^0 = \frac{1}{\delta x}\int_{x_i}^{x_{i+1}} w(0,x)dx,\ i \in \mathbb{Z} \qquad (3.3)$$

$$w_{i-1/2}^0 = \frac{1}{\delta x}\int_{x_{i-1}}^{x_i} w(0,x)dx,\ i \in \mathbb{Z} \qquad (3.4)$$

$$\rho_i^{k+1} = \rho_i^k + \frac{\delta t}{\delta x}\left(G_{i-1/2}^k - G_{i+1/2}^k\right), \; i \in \mathbb{Z},\, k \in \mathbb{Z},\, k \geq 0 \tag{3.5}$$

$$G_{i+1/2}^k = \frac{w_{i+1/2}^k + \left|w_{i+1/2}^k\right|}{2}\rho_i^k + \frac{w_{i+1/2}^k - \left|w_{i+1/2}^k\right|}{2}\rho_{i+1}^k,\; i \in \mathbb{Z},\, k \in \mathbb{Z},\, k \geq 0 \tag{3.6}$$

$$G_{i-1/2}^k = \frac{w_{i-1/2}^k + \left|w_{i-1/2}^k\right|}{2}\rho_{i-1}^k + \frac{w_{i-1/2}^k - \left|w_{i-1/2}^k\right|}{2}\rho_i^k,\; i \in \mathbb{Z},\, k \in \mathbb{Z},\, k \geq 0 \tag{3.7}$$

$$G_i^k = \frac{1}{2}\left(G_{i-1/2}^k + G_{i+1/2}^k\right),\; i \in \mathbb{Z},\, k \in \mathbb{Z},\, k \geq 0 \tag{3.8}$$

$$\rho_{i+1/2}^{k+1} = \frac{1}{2}\left(\rho_i^{k+1} + \rho_{i+1}^{k+1}\right),\; i \in \mathbb{Z},\, k \in \mathbb{Z},\, k \geq 0 \tag{3.9}$$

$$\hat{\beta}_i^k = \beta(\hat{w}_i^k) = \begin{cases} \beta(w_{i-1/2}^k), & G_i^k \geq 0 \\ \beta(w_{i+1/2}^k), & G_i^k < 0 \end{cases},\; i \in \mathbb{Z},\, k \in \mathbb{Z},\, k \geq 0 \tag{3.10}$$

$$\beta_{i+1/2}^k = \beta\left(w_{i+1/2}^k\right),\; i \in \mathbb{Z},\, k \in \mathbb{Z},\, k \geq 0 \tag{3.11}$$

$$\beta_{i+1/2}^{k+1} = \beta\left(w_{i+1/2}^{k+1}\right),\; i \in \mathbb{Z},\, k \in \mathbb{Z},\, k \geq 0 \tag{3.12}$$

$$\begin{aligned}\rho_{i+1/2}^{k+1}\beta_{i+1/2}^{k+1} &= \frac{1}{1+\sigma\delta t}\left(\frac{\rho_i^k + \rho_{i+1}^k}{2}\beta_{i+1/2}^k + \frac{\delta t}{\delta x}\left(P(\rho_i^{k+1}) - P(\rho_{i+1}^{k+1}) + G_i^k\hat{\beta}_i^k - G_{i+1}^k\hat{\beta}_{i+1}^k\right)\right. \\ &\left. + \frac{\delta t}{\delta x}\left(\rho_{i+1}^{k+1}\kappa(\rho_{i+1}^{k+1})\frac{w_{i+3/2}^k - w_{i+1/2}^k}{\delta x} - \rho_i^{k+1}\kappa(\rho_i^{k+1})\frac{w_{i+1/2}^k - w_{i-1/2}^k}{\delta x}\right)\right),\; i \in \mathbb{Z},\, k \in \mathbb{Z},\, k \geq 0\end{aligned} \tag{3.13}$$

with $\delta t, \delta x > 0$ being the time-step and spatial discretization step and the functions $\beta(w)$, $P(\rho)$ being given by (2.15), (2.26), respectively. Moreover, $\rho_i^k$, $\rho_{i+1/2}^k$, $w_i^k$, $w_{i+1/2}^k$ are the numerical values at the points $x_i = i\delta x$, $x_{i+1/2} = (i+1/2)\delta x$ and at time $t_k = k\delta t$ while $\rho_i^{k+1}$, $\rho_{i+1/2}^{k+1}$, $w_i^{k+1}$, $w_{i+1/2}^{k+1}$ are the numerical values at the points $x_i = i\delta x$, $x_{i+1/2} = (i+1/2)\delta x$, and at time $t_{k+1} = (k+1)\delta t$.

The next numerical method that we use, is the so-called "method of lines" (see [10], [29]). In this method, the states (density and mean speed) are approximated at fixed grid points in Lagrangian coordinates. Therefore, in order to use this numerical method, it is required to transform equations (2.1), (2.2) into Lagrangian coordinates by making the following transformation (see also [10]):

$$s = \int_{\tilde{l}(\tau)}^{\xi} \tilde{\rho}(\tau, z)dz,\; h = \tau. \tag{3.14}$$

Thus, in the Lagrangian coordinates, equations (2.1), (2.2) become

$$\tilde{\rho}_h + \tilde{\rho}^2\tilde{v}_s = 0 \tag{3.15}$$

$$(\tilde{\beta}(\tilde{v}))_h + (\tilde{P}(\tilde{\rho}))_s = (\tilde{\rho}^2\tilde{\kappa}(\tilde{\rho})\tilde{v}_s)_s - \tilde{\sigma}\tilde{\beta}(\tilde{v}) \tag{3.16}$$

for all $h > 0$ and $0 < s < \tilde{m}$ where $\tilde{m} = \int_{\tilde{l}(\tau)}^{\tilde{L}(\tau)} \tilde{\rho}(\tau, z) dz$, $\tilde{P}(\tilde{\rho}) = \tilde{\sigma} \int_{\bar{\rho}}^{\tilde{\rho}} \tilde{\kappa}(\tau) d\tau$ for all $\tilde{\rho} \in (0, \rho_{max})$ and $\tilde{\beta}(\tilde{v})$ defined by (2.5). Model (3.15), (3.16) is studied under the following boundary conditions:

$$\tilde{\rho}(h, 0) = \tilde{\rho}(h, \tilde{m}) = 0. \tag{3.17}$$

The "method of lines" is applied by discretizing the spatial derivatives of model (3.15), (3.16) with respect to $s$ and leaving the time variable continuous. Therefore, the following system of $2n$ ordinary differential equations is obtained:

$$\frac{d}{dh}\tilde{\rho}_{2i}(h) + \{\tilde{\rho}_{2i}(h)\}^2 \frac{\tilde{v}_{2i+1}(h) - \tilde{v}_{2i-1}(h)}{\delta s} = 0 \tag{3.18}$$

$$\frac{d}{dh}\tilde{\beta}(\tilde{v}_{2i}(h)) + \frac{\tilde{P}(\tilde{\rho}_{2i}(h)) - \tilde{P}(\tilde{\rho}_{2i-2}(h))}{\delta s} = \frac{1}{(\delta s)^2}\{G_{2i}(h)(\tilde{v}_{2i+1}(h) - \tilde{v}_{2i-1}(h)) - G_{2i-2}(h)(\tilde{v}_{2i-1}(h) - \tilde{v}_{2i-3}(h))\} - \tilde{\sigma}\tilde{\beta}(\tilde{v}_{2i-1}(h)) \tag{3.19}$$

with boundary conditions

$$\tilde{\rho}_0(h) = \tilde{\rho}_{2n}(h) = 0 \tag{3.20}$$

and initial data

$$\tilde{\rho}_{2i}(0) = \tilde{\rho}_0\left(2i\frac{\delta s}{2}\right), \quad \tilde{v}_{2i-1}(0) = \tilde{v}_0\left((2i-1)\frac{\delta s}{2}\right) \tag{3.21}$$

where $i = 1, 2, ..., n$, $G_{2i}(h) = \tilde{\rho}_{2i}^2(h)\tilde{\kappa}(\tilde{\rho}_{2i}(h))$, $\tilde{v}_{-1}(h) = \tilde{v}_{2n+1}(h) = v^*$ and $\delta s = \frac{\tilde{m}}{n}$ for a given number $n$ of grid-points.

## 4. Numerical Experiments

*4.1 Academic Example*

In this example, we apply the particle method (2.32)-(2.37), and use (2.48)-(2.54) to approximate the solutions of model (2.11)-(2.14). Moreover, we compare the numerical solutions obtained by applying the proposed particle method with the numerical solutions of the modified numerical methods presented in Section 3. Finally, a comparison with the reduced model (2.20) is included, indicating the effectiveness of the employed particle method as well as the high-resolution approximations provided by the reduced model.

In order to approximate the solutions of the ODE systems (2.32)-(2.37) and (3.18)-(3.21) involved in the particle method and the method of lines, respectively, we use an adaptive step size numerical scheme (see [12], pages 167-169) consisting of the explicit Euler and Heun schemes. For these schemes, the new step size is given by the formula:

$$\delta t_{new} = \delta t \min\left\{p, 0.9\sqrt{\frac{1}{err}}\right\} \tag{4.1}$$

where

$$err = \sqrt{\frac{1}{2n}\sum_{i=1}^{n}\left(\frac{x_{i,Euler} - x_{i,Heun}}{sc_{x,i}}\right)^2 + \frac{1}{2n}\sum_{i=1}^{n}\left(\frac{y_{i,Euler} - y_{i,Heun}}{sc_{y,i}}\right)^2} \quad (4.2)$$

and

$$sc_{x,i} = Atol + Rtol \max\{|x_i|, |x_{i,Heun}|\}, \quad i = 1,...,n \quad (4.3)$$

$$sc_{y,i} = Atol + Rtol \max\{|y_i|, |y_{i,Heun}|\}, \quad i = 1,...,n \quad (4.4)$$

for a given number $n$ of state variables.

The constant $Atol > 0$ denotes the tolerance for absolute errors, $Rtol > 0$ is the tolerance for relative errors, and $p \geq 1$ is a constant factor which determines the magnitude of a possible increase of the step size. In addition, $x_{i,Euler}$, $y_{i,Euler}$ and $x_{i,Heun}$, $y_{i,Heun}$, $i = 1,...,n$ are the components of the solution by the respective schemes, where $(x_i, y_i) = (x_i, w_i)$, $i = 1,...,n$ for the ODE system (2.32)-(2.37) and $(x_i, y_i) = (\tilde{\rho}_{2i}, \tilde{v}_{2i})$, $i = 1,...,n$ for the ODE system (3.18)-(3.21). For the particle method (2.32)-(2.37) we use the aforementioned adaptive numerical method with $A_{tol} = R_{tol} = 10^{-4}$, $p = 2$ and the parameters of Table 3. To obtain the functions $\Phi'(s)$ and $K(s)$, we used (2.38), (2.42) and the function

$$\kappa(\rho) = \begin{cases} 0, 0 < \rho \leq 1 \\ \dfrac{(\rho-1)^2}{\rho(R-\rho)}, 1 < \rho < R \end{cases} \quad (4.5)$$

for the state-dependent viscosity term in (2.11)-(2.14). Finally, to obtain the approximations of the solutions of (2.11)-(2.14) we use the equations (2.48)-(2.54) and the initial condition:

$$\rho[0] = \begin{cases} 0.25(x-l(0))^2(x-L(0))^2, & l(0) < x < L(0) \\ 0, & otherwise \end{cases}$$

$$w[0] = \begin{cases} -0.158(x-0.5)^2(x-1.5)^2, & 0.5 < x < 1.5 \\ 0, & x \in (l(0), 0.5] \cup [1.5, L(0)) \end{cases} \quad (4.6)$$

where $l(0) = -0.52$ and $L(0) = 2.52$.

|  | $n$ | $b$ | $R$ | $\sigma$ | $a$ |
|---|---|---|---|---|---|
| Particle method (2.32)-(2.37) | 225 | 0.0606 | 1.9 | 30 | 0.4653 |

**Table 3**: Parameters of the Particle Method (2.32)-(2.37).

Figure 1, shows the solution profiles and their convergence to an equilibrium point $(\rho, w)$ with $0 \leq \rho(x) \leq 1$ and $w(x) = 0$ for all $x \in (l(t), L(t))$ as time increases. Figure 2 and Figure 3 illustrate the time evolution of the discretized functionals defined by (2.43), (2.44). Notice that both functionals $E_n$, $W_n$ are not increasing, showing that the differential inequalities (2.46), (2.47) are preserved.

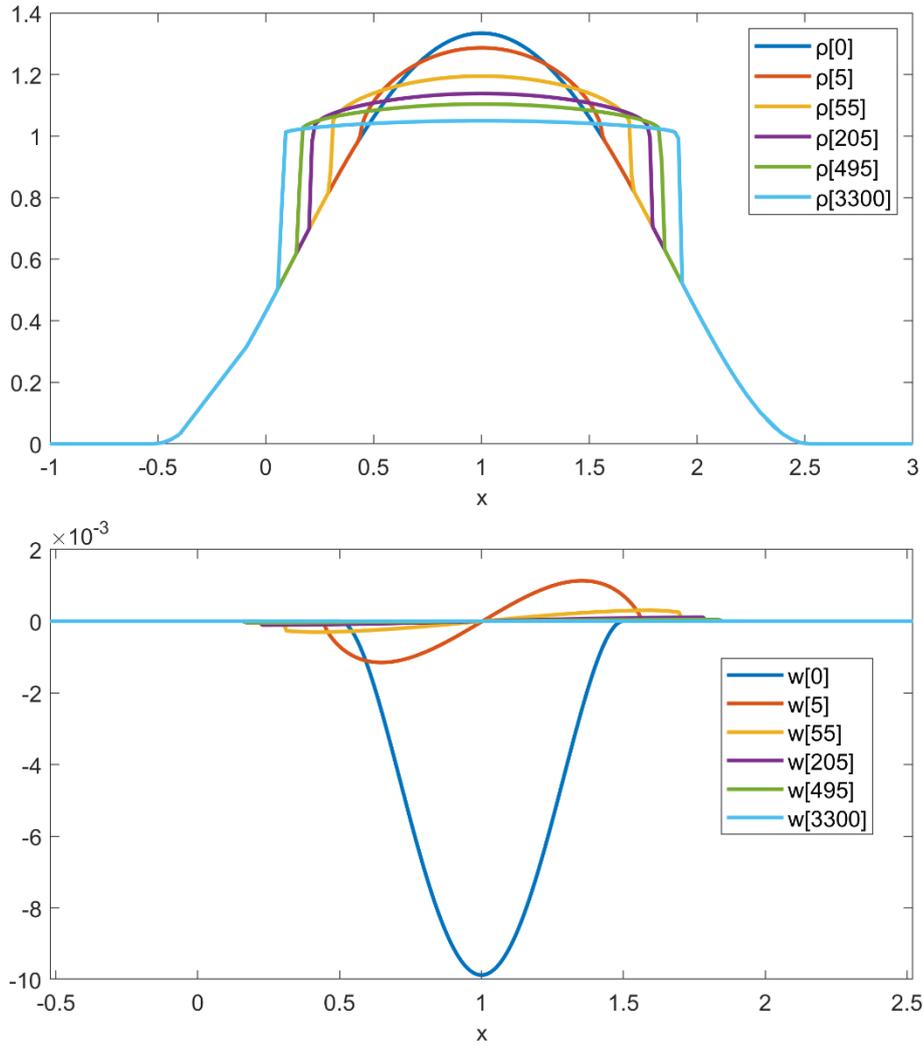

**Figure 1**: Solution profiles of model (2.11)-(2.14) produced by applying the particle method (2.32)-(2.37).

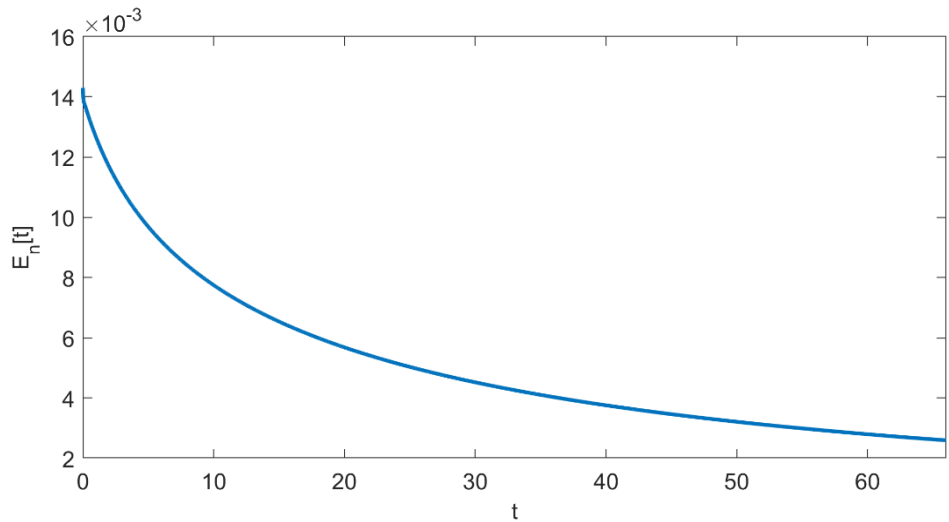

**Figure 2**: Time-evolution of the discretized functional $E_n$ defined by (2.43) obtained by using the particle method (2.32)-(2.37).

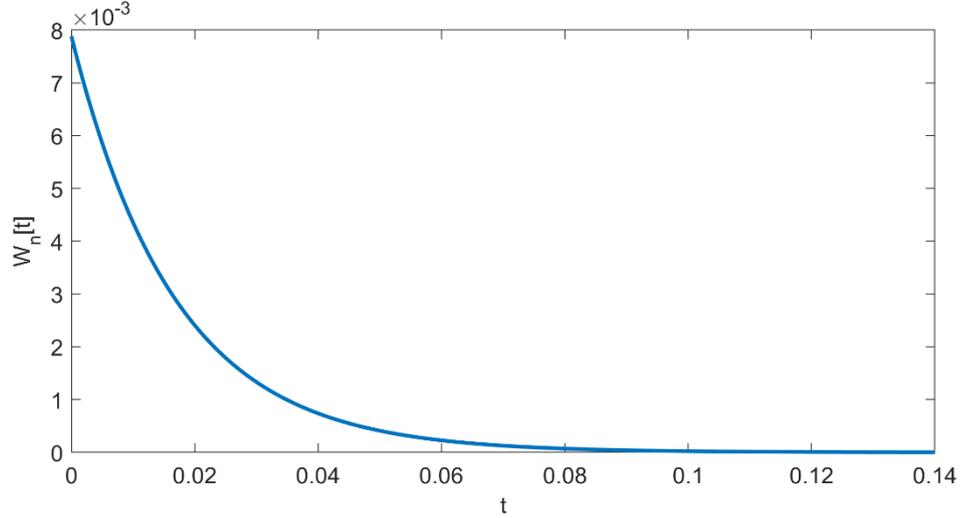

**Figure 3**: Time-evolution of the discretized functional $W_n$ defined by (2.44) obtained by using the particle method (2.32)-(2.37).

Next, we compare the numerical solutions produced by applying the proposed particle method with those of the numerical method (3.2)-(3.13) (NM1) and the method of lines (3.18)-(3.21) (NM2). For the numerical method (3.2)-(3.13) (NM1), we used the parameters of Table 4. To compare the solutions of the particle method with the solutions of the method of lines (3.18)-(3.21) (NM2), we transform the solution of (3.18)-(3.21) from Lagrangian to Eulerian coordinates and remove its physical dimensions through the variable transformation $\xi = rx + v^*\tau$ and the dimensionless quantities $\tau = \frac{r}{v^*}t$, $b = \frac{v_{max} - v^*}{v^*}$, $R = \frac{\rho_{max}}{\bar{\rho}} > 1$, $w = \frac{\tilde{v} - v^*}{v^*}$, $\rho = \frac{\tilde{\rho}}{\bar{\rho}}$, $\sigma = \frac{r}{v^*}\tilde{\sigma}$, where the values of $\tilde{\sigma}$, $v^*$, $r$, $\bar{\rho}$, $\rho_{max}$, $v_{max}$ are given in Table 5. Note that the transformation was made by integrating the differential equation $\dot{\xi}_i = \tilde{v}_{2i}$, $i = 1,...,n$, using the trapezoidal rule.

To obtain the numerical solution of (3.18)-(3.21) we use again the aforementioned adaptive step-size scheme with $A_{tol} = R_{tol} = 10^{-4}$, $p = 2$. Note also that the expressions of the functions $\tilde{\beta}(\tilde{v})$, $\tilde{P}(\tilde{\rho})$, $\tilde{\kappa}(\tilde{\rho})$ used in the ODE system (3.18)-(3.21) can be found by using the function $\tilde{P}(\tilde{\rho}) = \tilde{\sigma}\int_{\bar{\rho}}^{\tilde{\rho}} \tilde{\kappa}(\tau)d\tau$, for all $\tilde{\rho} \in (0, \rho_{max})$, the definitions (4.5), (2.15) and the relations of Table 1.

| | $n$ | $\delta t$ | $\delta x$ | $\sigma$ |
|---|---|---|---|---|
| Numerical method 1 (NM1) (3.2)-(3.13) | 225 | 0.0033 | 0.0178 | 30 |

**Table 4**: Parameters of the numerical method (3.2)-(3.13) (NM1).

| | $n$ | $\tilde{\sigma}$ | $\delta s$ | $v^*$ | $r$ | $\bar{\rho}$ | $\rho_{max}$ | $v_{max}$ |
|---|---|---|---|---|---|---|---|---|
| Numerical Method 2 (NM2) (Method of lines) (3.18)-(3.21) | 225 | 990 | 1.3665 | 33 | 1 | 63.1579 | 120 | 35 |

**Table 5**: Parameters of the numerical method (3.18)-(3.21) (NM2).

Figures 4 and 5 illustrate the small difference between the discretized functionals $E_n[t]$, $W_n[t]$ obtained by applying the proposed particle method (2.32)-(2.37) and the numerical methods (3.2)-(3.13) (NM1), (3.18)-(3.21) (NM2). Moreover, the sup-norms of the difference between the numerical solutions for the states $\rho, w$ are depicted in Figures 6, 7. It is obvious that as time increases, and the solution converges to the equilibrium point where $0 \leq \rho(x) \leq 1$, $w(x) = 0$, the difference between the numerical solutions becomes even smaller.

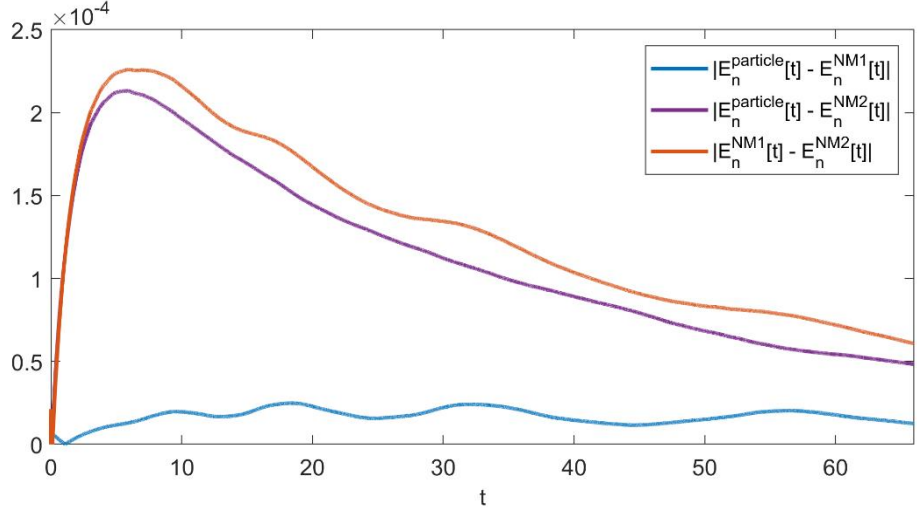

**Figure 4**: Time-evolution of $\left|E_n^{particle}[t] - E_n^{NM1}[t]\right|$, $\left|E_n^{particle}[t] - E_n^{NM2}[t]\right|$, $\left|E_n^{NM1}[t] - E_n^{NM2}[t]\right|$.

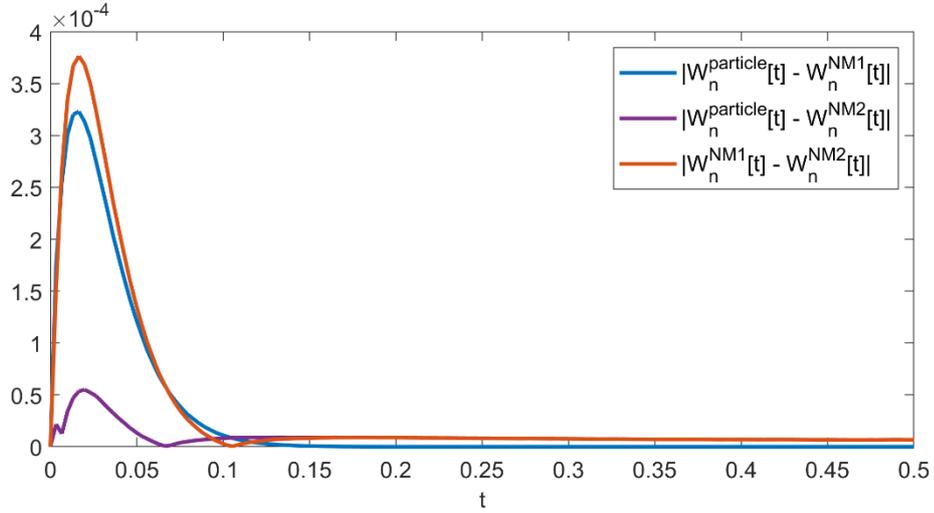

**Figure 5**: Time-evolution of $\left|W_n^{particle}[t] - W_n^{NM1}[t]\right|$, $\left|W_n^{particle}[t] - W_n^{NM2}[t]\right|$ $\left|W_n^{NM1}[t] - W_n^{NM2}[t]\right|$.

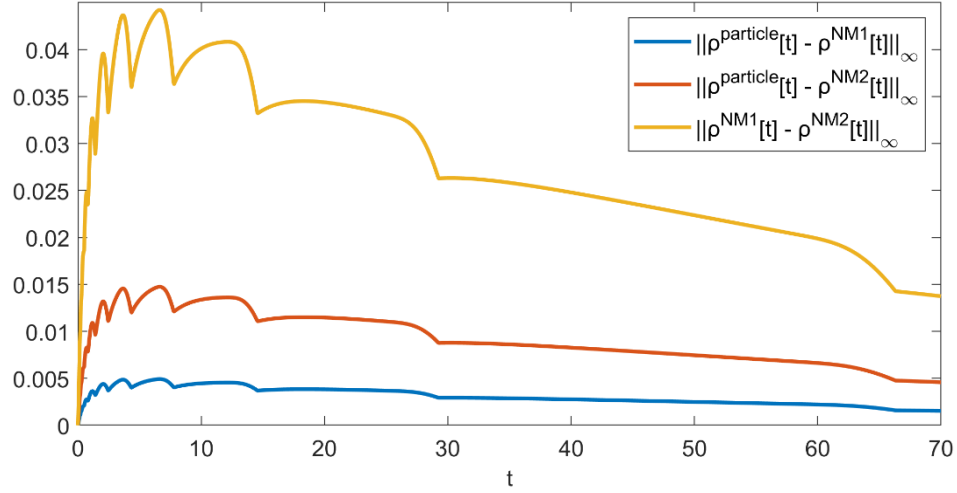

**Figure 6**: Time-evolution of $\left\|\rho^{particle}[t]-\rho^{NM1}[t]\right\|_\infty$, $\left\|\rho^{particle}[t]-\rho^{NM2}[t]\right\|_\infty$, $\left\|\rho^{NM1}[t]-\rho^{NM2}[t]\right\|_\infty$.

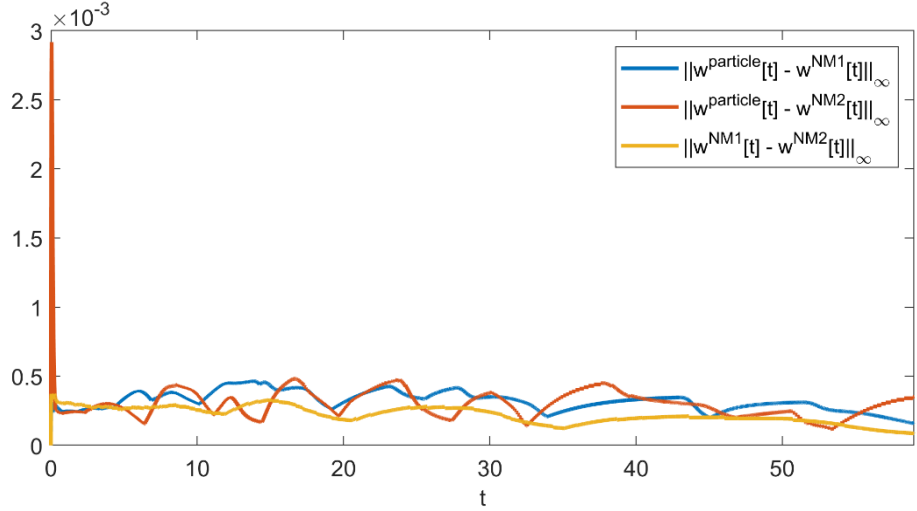

**Figure 7**: Time-evolution of $\left\|w^{particle}[t]-w^{NM1}[t]\right\|_\infty$, $\left\|w^{particle}[t]-w^{NM2}[t]\right\|_\infty$, $\left\|w^{NM1}[t]-w^{NM2}[t]\right\|_\infty$.

Notice that (2.47) along with (2.44) imply the differential equation $\dot{W}_n(t) = -2\sigma W_n(t)$ with analytical solution $W_n(t) = W_n(0)\exp(-2\sigma t)$. Figure 8 shows that the discretized energy functional $W_n$ obtained by using the particle method is closer to the analytical solution of the differential equation $\dot{W}_n(t) = -2\sigma W_n(t)$. Indeed, by plotting the time evolution of $\ln(W_n(t)/W_n(0))$ in Figure 9, it is shown that when the values of $W_n$ become very small, the discretized functional $W_n$ produced by applying the particle method (2.32)-(2.37) remains closer to the analytical solution $W_n(t) = W_n(0)\exp(-2\sigma t)$ compared to the discretized functionals of the numerical methods (3.2)-(3.13) (NM1), (3.18)-(3.21) (NM2), which diverge. Note also that a better approximation of the solution can be achieved with the particle method by using smaller values for the tolerances $Atol, Rtol$.

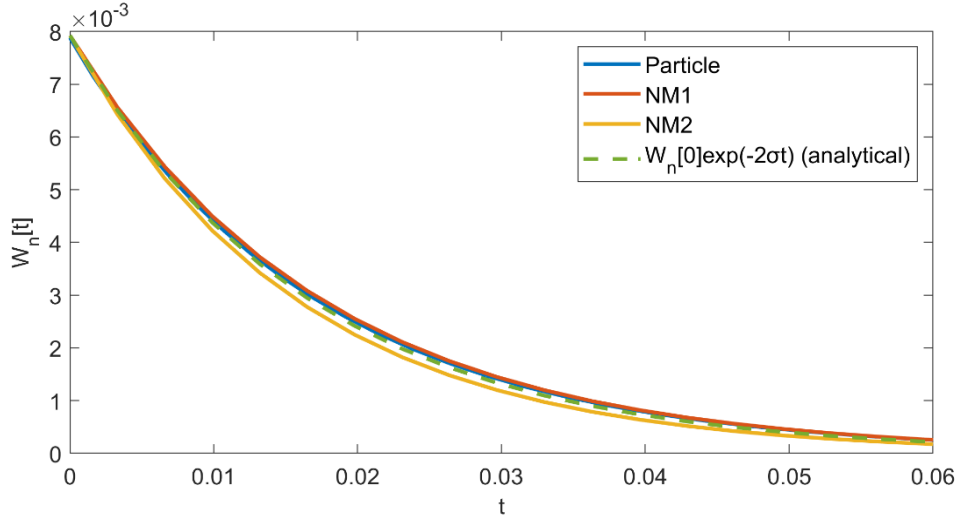

**Figure 8**: Time-evolution of $W_n^{particle}[t]$, $W_n^{NM1}[t]$, and $W_n^{NM2}[t]$.

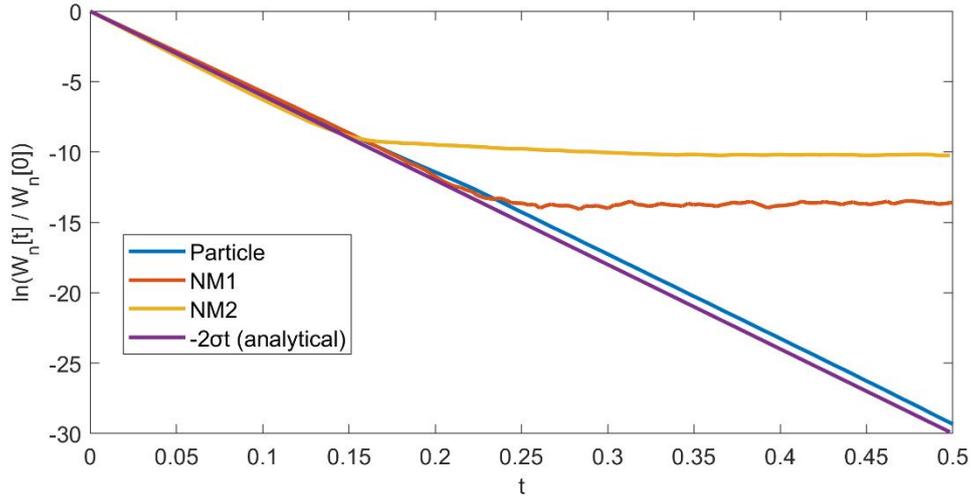

**Figure 9**: Time-evolution of $\ln(W_n^{particle}[t]/W_n[0])$, $\ln(W_n^{NM1}[t]/W_n[0])$, and $\ln(W_n^{NM2}[t]/W_n[0])$.

In addition, we compare the computational cost for each numerical method by means of CPU time needed. Table 6 shows the corresponding results for the overall simulation time $t_{max}$ of 33 time units, demonstrating that the CPU time needed for the proposed particle method is significantly less than the corresponding CPU times of the numerical methods (3.2)-(3.13) (NM1), (3.18)-(3.21) (NM2). It should be mentioned, that even in the case of smaller tolerances ($A_{tol} = R_{tol} = 10^{-8}$) for the particle method, its execution time was no more than 20 seconds. On the other hand, NM2 failed to give a result even after 2h. The simulations were performed on Matlab R2021a and a laptop with Intel(R) Core(TM) i7-7700HQ CPU 2.80GHz.

| Numerical Method | CPU time (sec) |
|---|---|
| Particle Method (2.32)-(2.37) | 18.65 |
| NM1: (3.2)-(3.13) | 274.52 |
| NM2: (3.18)-(3.21) | 171.97 |

**Table 6**: CPU times needed for the proposed particle method (2.32)-(2.37) and the numerical methods (3.2)-(3.13), (3.18)-(3.21) with $t_{max} = 33$ time units.

Finally, we compare the numerical solution of model (2.11)-(2.14) produced by applying the particle method (2.32)-(2.37) with the numerical solution of the reduced model (2.20) obtained by using the numerical scheme proposed in [34] with time-step $\delta t = 0.0033$ and spatial discretization step $\delta x = 0.0178$. Figure 10 and Figure 11 show that the selection of the constant $\sigma$ plays an important role in the deviation between the numerical solutions of the two models.

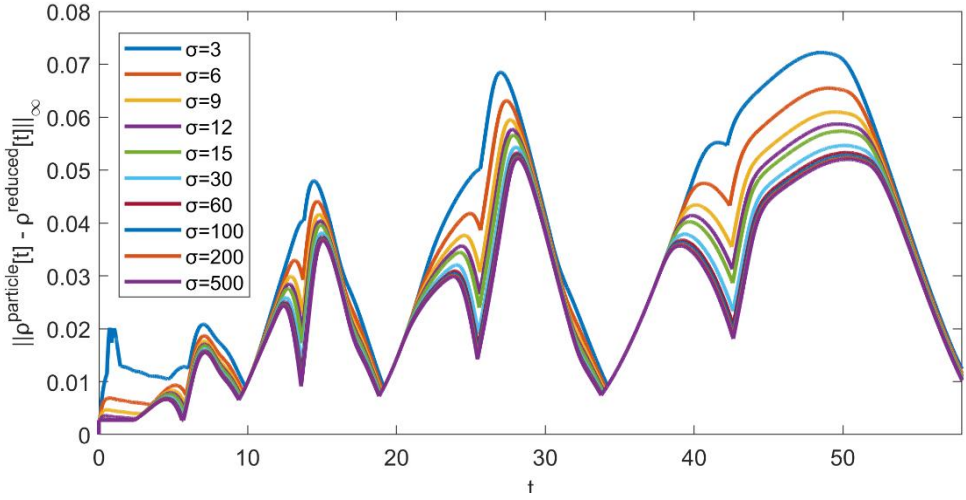

**Figure 10**: Time-evolution of $\left\|\rho^{particle}[t] - \rho^{reduced}[t]\right\|_\infty$ for a range of $\sigma$ values.

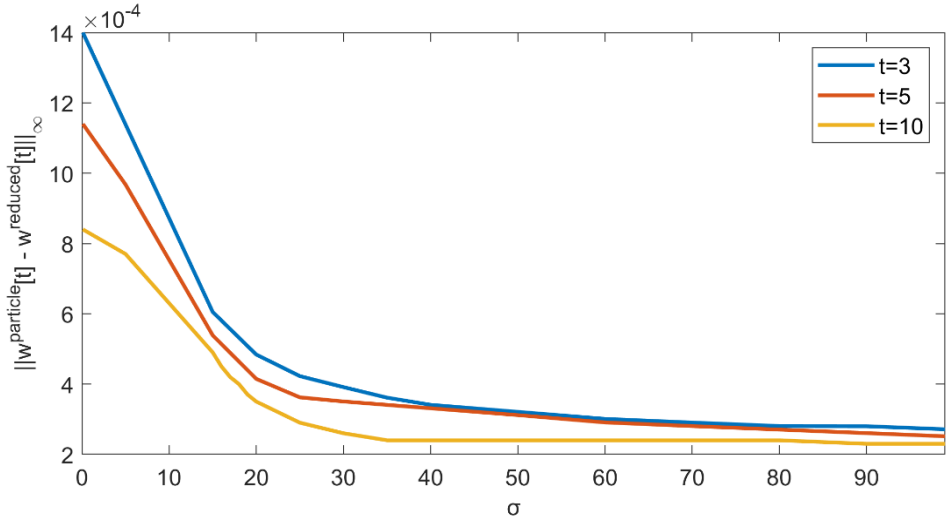

**Figure 11**: Profiles of $\left\|w^{particle}[t] - w^{reduced}[t]\right\|_\infty$ for a range of $\sigma$ values.

As the value of constant $\sigma$ increases, the convergence of the solution of model (2.11)-(2.14) to the solution of the reduced model becomes faster and therefore the difference between the numerical solutions is faster suppressed. Figure 12 illustrates the small difference between the numerical solution of model (2.11)-(2.14) and the numerical solution of the reduced model (2.20). Although a different initial condition $w[0]$ is given for each model, it is obvious that the sup-norm of the difference between the numerical solutions becomes small very fast, indicating the high-resolution approximation provided by the reduced model.

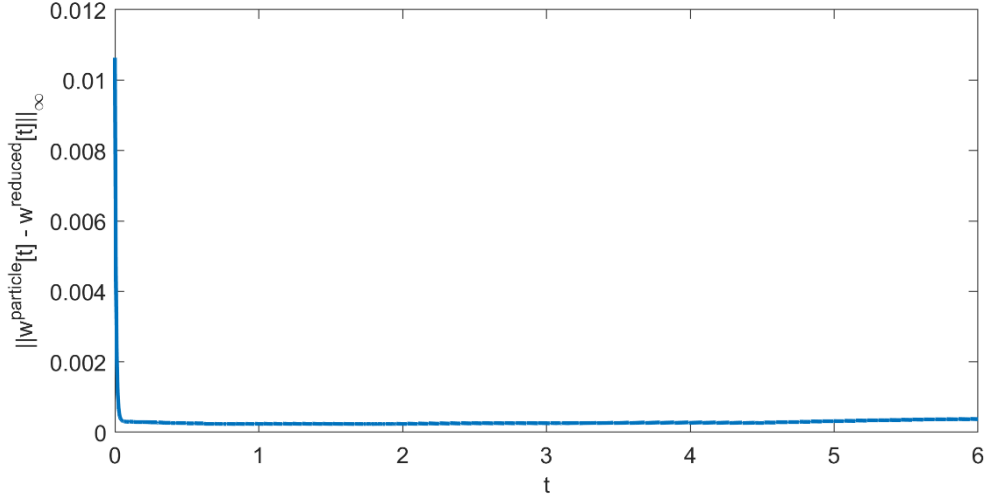

**Figure 12**: Time evolution of $\left\| w^{particle}[t] - w^{reduced}[t] \right\|_{\infty}$.

*4.2 Traffic Application*

In this section, a traffic scenario is considered to compare the density, mean speed, acceleration, and flow of the macroscopic model for automated vehicles (AV model) (2.1)-(2.4) with a well-known second-order macroscopic traffic flow model for human drivers, the so-called ARZ (Aw-Rascle-Zhang) model (see [3], [35]). For this scenario, we also compare the numerical solutions produced by applying the particle method (2.32)-(2.37) with the numerical solutions of (3.2)-(3.13) (NM1), (3.18)-(3.21) (NM2) and the numerical solutions of the reduced model (2.10). Note that this example is an analogue of the traffic application presented in [34], where the reduced AV model (2.10) was compared with a first-order traffic flow model for human drivers, the so-called LWR (Lighthill-Witham-Richards) model (see [25], [30]). The form of the ARZ model is given as follows:

$$\tilde{\rho}_\tau + (y + \tilde{\rho} V(\tilde{\rho}))_\xi = 0 \tag{4.7}$$

$$y_\tau + (y^2/\tilde{\rho} + y V(\tilde{\rho}))_\xi = -\frac{y}{\delta} \tag{4.8}$$

where $y = \tilde{\rho}(\tilde{v} - V(\tilde{\rho}))$ and $V(\tilde{\rho})$ is the equilibrium speed profile given by the expression:

$$V(\tilde{\rho}) = v_f \exp\left(-\frac{1}{d}\left(\frac{\tilde{\rho}}{\rho_c}\right)^d\right). \tag{4.9}$$

Note that $d$ is a positive constant, $\rho_c$ is the critical density, $v_f$ is the free speed and $\delta$ is the relaxation time. In this scenario, a single-lane motorway is considered with $v_f = 102\,(km/h)$, $\rho_{max} = 180\,(veh/km)$, $d = 2.34$ and $\rho_c = 33.3\,(veh/km)$ that were estimated in [23] based on real data from a part of the Amsterdam A10 motorway. The initial condition used for this scenario is illustrated in Figure 13. Notice that the initial density on the road is characterized by a congestion belt in the interval $[1.5, 2.75]\,(km)$ with $\tilde{\rho}_0(\xi) = 0$ for $\xi \leq \tilde{l}(0) = 0$ and $\xi \geq \tilde{L}(0) = 4$. Moreover, the density has been selected to be lower than the critical density $\rho_c$ outside the congested area.

For the AV model (2.1)-(2.4), we apply the particle method (2.32)-(2.37) by using the adaptive step size numerical scheme discussed in Section 4.1 with $p = 2$ and tolerances $Atol = Rtol = 10^{-5}$. We used the following function:

$$\kappa(\rho) = c \begin{cases} 0, 0 < \rho \leq 1 \\ \dfrac{(\rho-1)^2}{\rho(R-\rho)}, 1 < \rho < R \end{cases} \qquad (4.10)$$

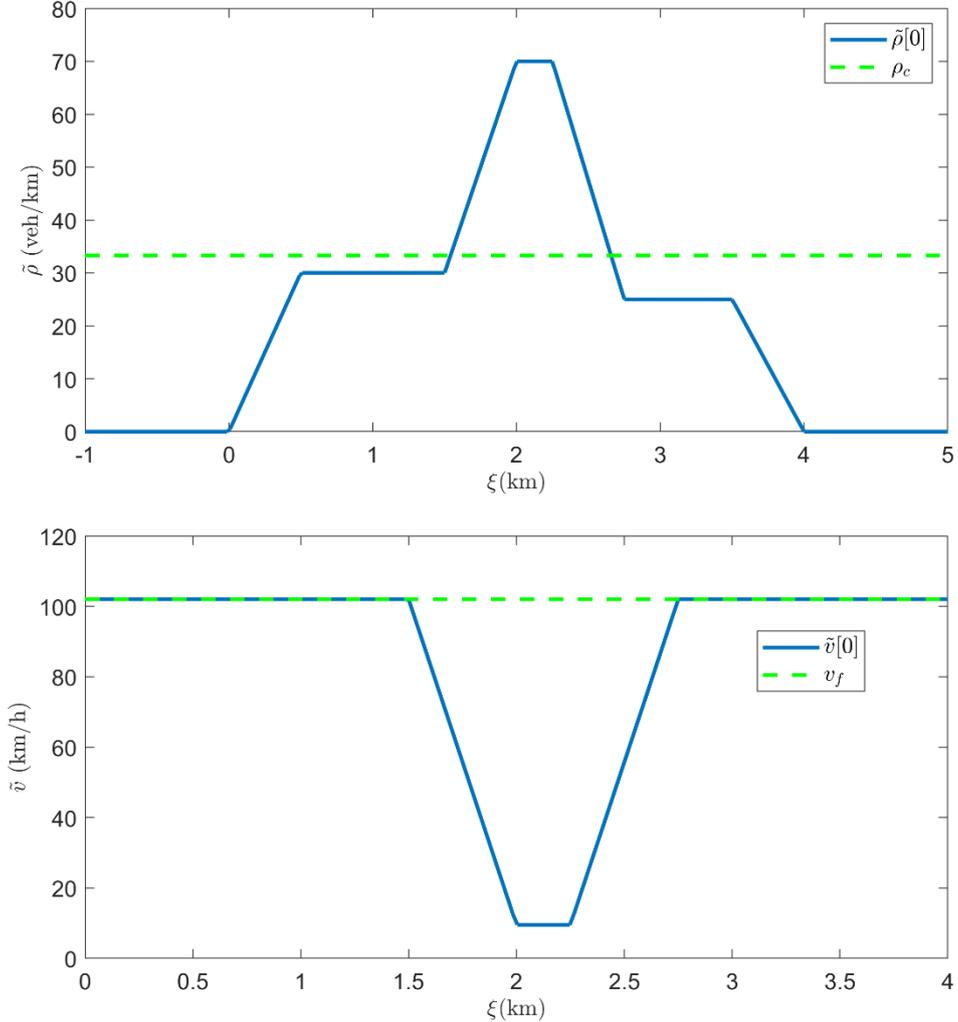

**Figure 13**: Initial density $\tilde{\rho}_0(\xi)$ (top) and speed $\tilde{v}_0(\xi)$ (bottom).

where $c > 0$ is a constant that is used to increase or decrease the viscosity. In this example, we selected $c = 40$. Note that the functions $K(s), \Phi'(s)$ are given by means of (4.10) and (2.38), (2.42), while the expressions of the functions $\tilde{\beta}(\tilde{v})$, $\tilde{\kappa}(\tilde{\rho})$ involved in the AV model (2.1)-(2.4) can be found by using (4.10), (2.15) and the relations of Table 1. Moreover, we use the parameters $v^* = 102\,(km/h)$ (equal to the free speed $v_f$), $\bar{\rho} = 31\,(veh/km)$, $r = 1\,(km)$, $\tilde{\sigma} = 3060\,(1/h)$, and recall that the quantities $R$, $b$, and $\sigma$ are given by $R = \frac{\rho_{max}}{\bar{\rho}} > 1$, $b = \frac{v_{max} - v^*}{v^*}$, and $\sigma = \frac{r}{v^*}\tilde{\sigma}$, respectively.

To approximate the solutions of the ARZ model (4.7)-(4.9), we use a second-order central scheme for balance laws described in [27] which is the extension of the Nessyahu-Tadmor scheme for conservation laws proposed in [28]. By expressing the ARZ model (4.7)-(4.9) in the following conservative form:

$$U_\tau + (f(U))_\xi = R(U) \tag{4.11}$$

where $U = \begin{bmatrix} \tilde{\rho} \\ y \end{bmatrix}$, $f(U) = \begin{bmatrix} y + \tilde{\rho}V(\tilde{\rho}) \\ y^2/\tilde{\rho} + yV(\tilde{\rho}) \end{bmatrix}$, and $R(U) = \begin{bmatrix} 0 \\ -y/\delta \end{bmatrix}$, the scheme reads as follows:

$$U_{i+1/2}^{k+1} = \frac{1}{2}(U_i^k + U_{i+1}^k) + \frac{1}{8}(U_i' - U_{i+1}') - \frac{\delta\tau}{\delta\xi}(f(U_{i+1}^{k+1/2}) - f(U_i^{k+1/2}))$$
$$+ \delta\tau\left(\frac{3}{8}R(U_i^{k+1/3}) + \frac{3}{8}R(U_{i+1}^{k+1/3}) + \frac{1}{4}R(U_{i+1/2}^{k+1})\right) \tag{4.12}$$

where

$$U_i' = \min\mathrm{mod}(U_{i+1}^k - U_i^k, U_i^k - U_{i-1}^k) \tag{4.13}$$

$$U_i^{k+1/2} = U_i + \frac{\delta\tau}{2}\left(R(U_i^{k+1/2}) - \frac{f_i'}{\delta\xi}\right) \tag{4.14}$$

$$U_i^{k+1/3} = U_i + \frac{\delta\tau}{3}\left(R(U_i^{k+1/3}) - \frac{f_i'}{\delta\xi}\right) \tag{4.15}$$

$$f_i' = \min\mathrm{mod}(f(U_{i+1}) - f(U_i), f(U_i) - f(U_{i-1})) \tag{4.16}$$

and

$$\min\mathrm{mod}(x, y) = \begin{cases} \mathrm{sgn}(x)\min(|x|, |y|), & \text{if } \mathrm{sgn}(x) = \mathrm{sgn}(y) \\ 0, & \text{otherwise} \end{cases} \tag{4.17}$$

By selecting $\delta = 20\,(s)$ (the same value is used in [3]), $\delta\xi = 0.04$, $\delta\tau = 10^{-6}$ for the numerical scheme (4.12)-(4.17) and $a = 0.2411$, $n = 4500$ particles/vehicles on the interval $[0, 120]\,(km)$ for the particle method, we obtain the density and speed profiles illustrated in Figure 14 and Figure 15, respectively. Figure 14 shows that the dissipation of density is much stronger for the ARZ model while the density is spread over a large road interval for increasing $\tau > 0$. More specifically, for $\tau = 1\,(h)$, the density is non-zero over the interval $(90, 106)\,(km)$, which implies that the vehicles retain large inter-vehicle distances. Notice also that the density for the AV model converges towards an equilibrium where $\tilde{v}(\xi) = v^*$, $\tilde{\rho}(\xi) \leq \bar{\rho} < \rho_{cr}$ and the travel time is higher using the AV model, since, for $\tau = 1\,(h)$, vehicles are in the interval $(102, 106)\,(km)$, while with the ARZ model, vehicles are in the interval $(90, 106)\,(km)$.

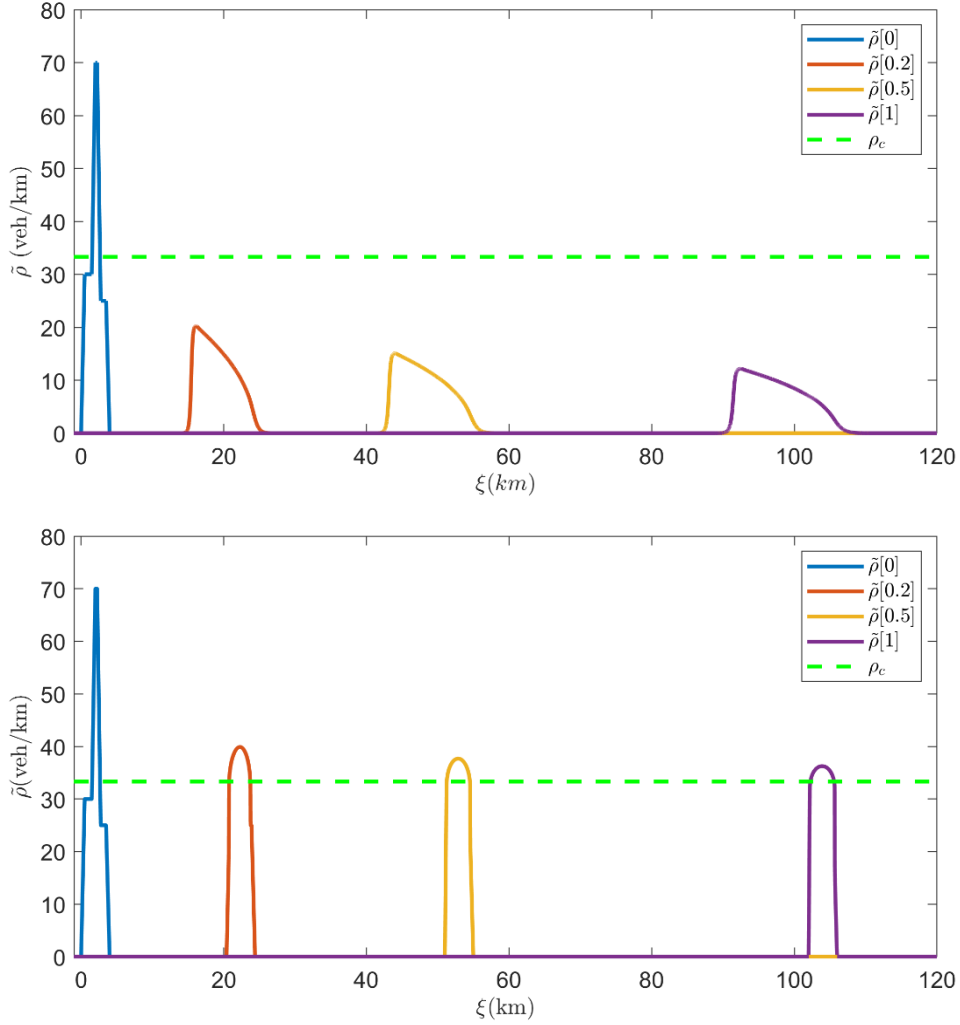

**Figure 14**: Density profiles for the ARZ model obtained by using $\delta = 20\,(s)$ (top), and density profiles for the AV model (2.1)-(2.4) with $v^* = 102\,(km/h)$ (bottom).

Due to the fact that the use of the AV model (2.1)-(2.4) causes the convergence of density towards an equilibrium near the critical density, while vehicles are spread in a small space interval (see Figure 14), it is expected that the flow $q = \tilde{\rho}\tilde{v}$ is much higher for the AV model compared to the ARZ. Indeed, this can be noticed in Table 7 where we compare the mean flow for each model using the definition:

$$\text{Mean Flow} = \frac{1}{T}\int_0^T \int_{a(t)}^{b(t)} \frac{\tilde{\rho}(t,x)\tilde{v}(t,x)}{b(t) - a(t)}\,dxdt \qquad (4.18)$$

where $T > 0$ denotes the time horizon and $[a(\tau), b(\tau)]$ denotes the interval where $\tilde{\rho}(\tau, \xi) \neq 0$ for $\xi \in [a(\tau), b(\tau)]$ at each time instant $\tau \geq 0$. Notice that the mean flow obtained with the ARZ model is achieved with the AV model (2.1)-(2.4) by using a significantly smaller desired speed $v^* = 30.8\,(km/h)$ than the free speed $v_f = 102\,(km/h)$ of the ARZ.

| | AV model (2.1)-(2.4) with $v^* = 102\,(km/h)$ | AV model (2.1)-(2.4) with $v^* = 70\,(km/h)$ | AV model (2.1)-(2.4) with $v^* = 30.8\,(km/h)$ | ARZ using $\delta = 20\,(s)$ $v_f = 102\,(km/h)$ |
|---|---|---|---|---|
| Mean Flow | 3388 (veh/h) | 2229 (veh/h) | 963 (veh/h) | 963 (veh/h) |

**Table 7**: Mean Flow

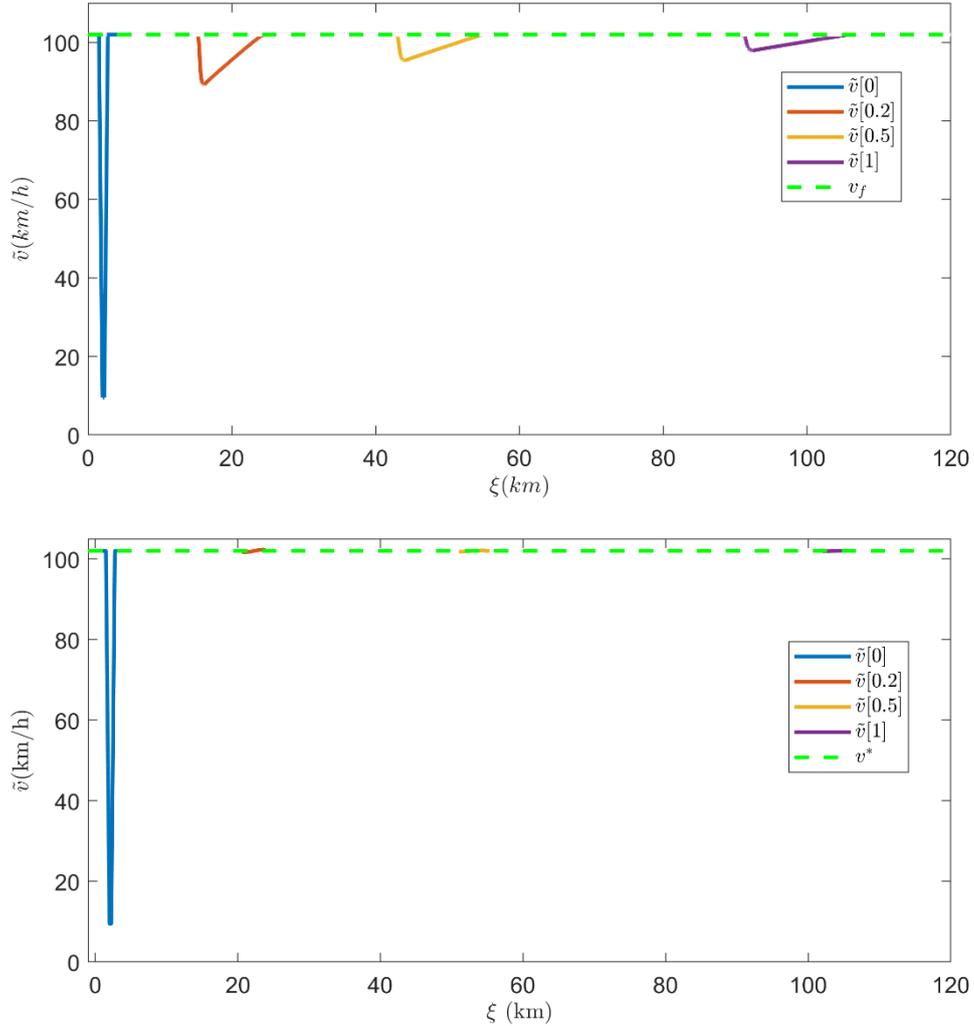

**Figure 15**: Speed profiles for the ARZ model obtained by using $\delta = 20\,(s)$ (top), and speed profiles of the AV model (2.1)-(2.4) with $v^* = 102\,(km/h)$ (bottom).

Figure 16 displays the sup-norms of the accelerations by using the ARZ and the AV model (2.1)-(2.4). Since the vehicles with the AV model obtain the desired speed $v^*$ faster compared to the ARZ (see Figure 15), the maximum accelerations of the vehicles using the AV model (2.1)-(2.4) are decreasing faster compared to those of the ARZ model, which fluctuate at the first time instants and then converge to lower values. Notice also that the acceleration of the AV model (2.1)-(2.4) is close to the acceleration of the ARZ model. This occurs due to the selection of the parameter $\tilde{\sigma}$ which appears on the friction term $-\tilde{\sigma}\tilde{\rho}\tilde{\beta}(\tilde{v})$ of the speed PDE (2.2). The constant $\tilde{\sigma}$ has the potential to affect the tendency of vehicles to adjust their speed to the speed set-point $v^*$. Consequently, by changing the value of the parameter $\tilde{\sigma}$ we can increase or decrease the acceleration of the vehicles.

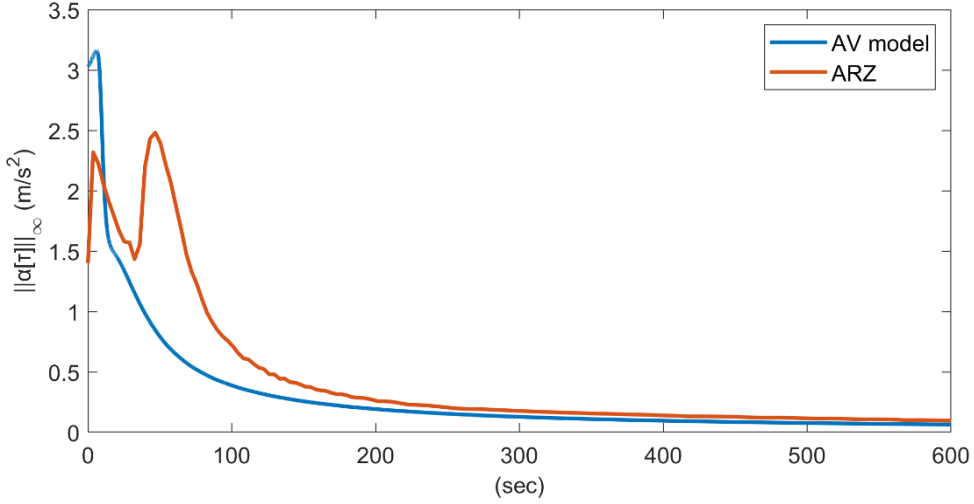

**Figure 16**: Sup-norm of the acceleration of the ARZ model and the AV model.

Next, we compare the numerical solutions produced by applying the particle method (2.32)-(2.37) with the numerical solutions of (3.2)-(3.13) (NM1), (3.18)-(3.21) (NM2) and the numerical solution of the reduced model (2.10). Figures 17, 18 illustrate the profile solutions of the AV model (2.1)-(2.4) produced by applying each one of the aforementioned numerical methods. It is obvious that the difference between the numerical solutions is very small and it tends to be even smaller as time increases and the solutions converge to the equilibrium point $\tilde{v}(\xi) \equiv v^*$, $\tilde{\rho}(\xi) \leq \bar{\rho}$, $\xi \in (\tilde{l}(\tau), \tilde{L}(\tau))$.

The computational cost of each numerical method, in terms of CPU time needed, is given in Table 8 for the overall simulation time of $t_{\max} = 0.1\,(h)$. Again, the CPU time needed for the proposed particle method is significantly less than the corresponding CPU times of NM1 (3.2)-(3.13), NM2 (3.18)-(3.21) and the numerical method proposed in [34] for the numerical approximation of the reduced model (2.10). The simulations were performed on Matlab R2021a and a laptop with Intel(R) Core(TM) i7-7700HQ CPU 2.80GHz.

| **Numerical Method** | **CPU time (sec)** |
|---|---|
| Particle Method (2.32)-(2.37) | 391.32 |
| NM1: (3.2)-(3.13) | 2034.76 |
| NM2: (3.18)-(3.21) | 1602.21 |
| Reduced model | 2104.13 |

**Table 8**: CPU times needed for the proposed particle method (2.32)-(2.37), the numerical methods (3.2)-(3.13) (NM1), (3.18)-(3.21) (NM2) and the numerical method proposed in [34] for the reduced model (2.10).

However, the proposed particle method (2.32)-(2.37) provides higher-resolution approximations compared to the numerical methods (3.2)-(3.13) (NM1), (3.18)-(3.21) (NM2). This conclusion can be reached by plotting the time evolution of $\ln(W_n(t)/W_n(0))$ in Figure 19. It is shown that the discretized functional $W_n$ obtained by using the particle method coincides with the analytical solution of the differential equation $\dot{W}_n(t) = -2\sigma W_n(t)$ while the discretized functionals produced by applying the methods (3.2)-(3.13) (NM1), (3.18)-(3.21) (NM2) diverge.

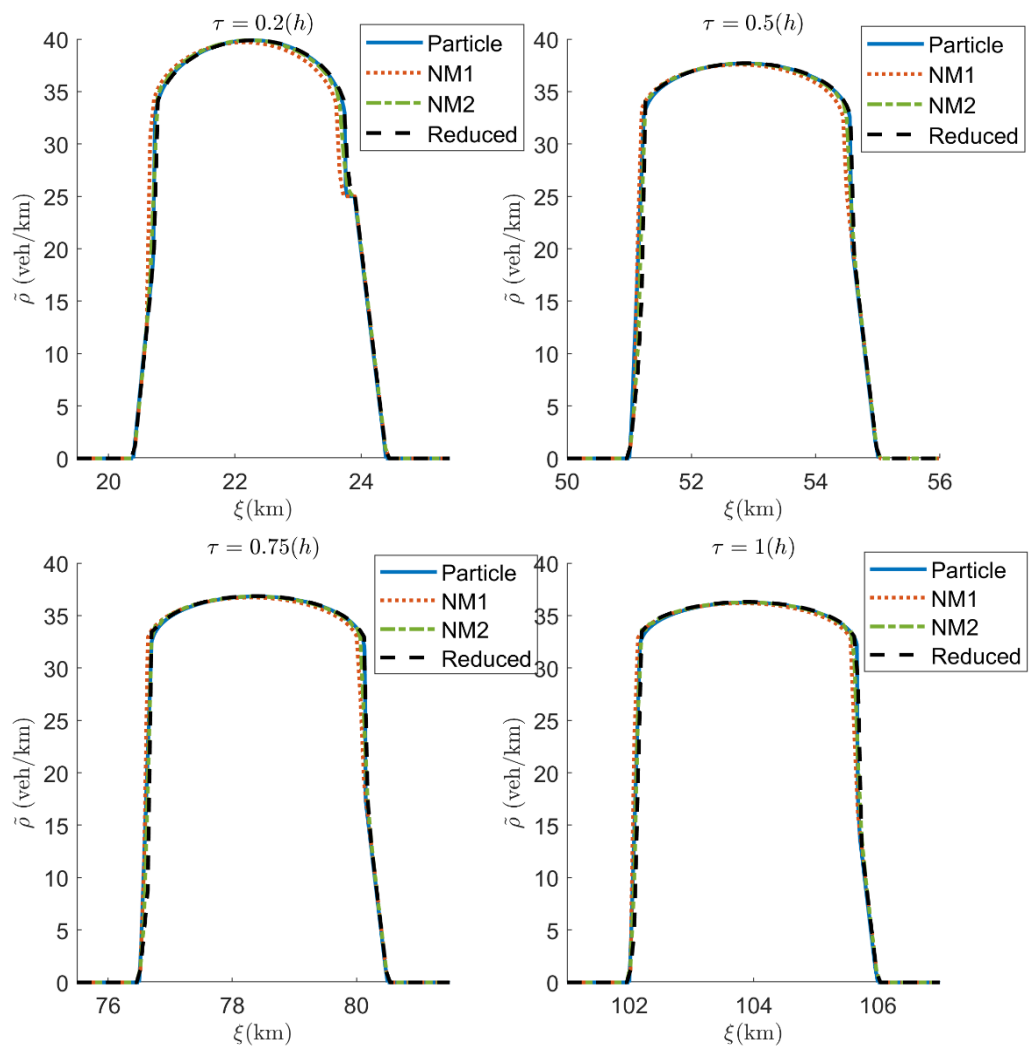

**Figure 17**: Density profiles of the AV model (2.1)-(2.4) at $\tau = 0.2(h)$, $\tau = 0.5(h)$, $\tau = 0.75(h)$, and $\tau = 1(h)$.

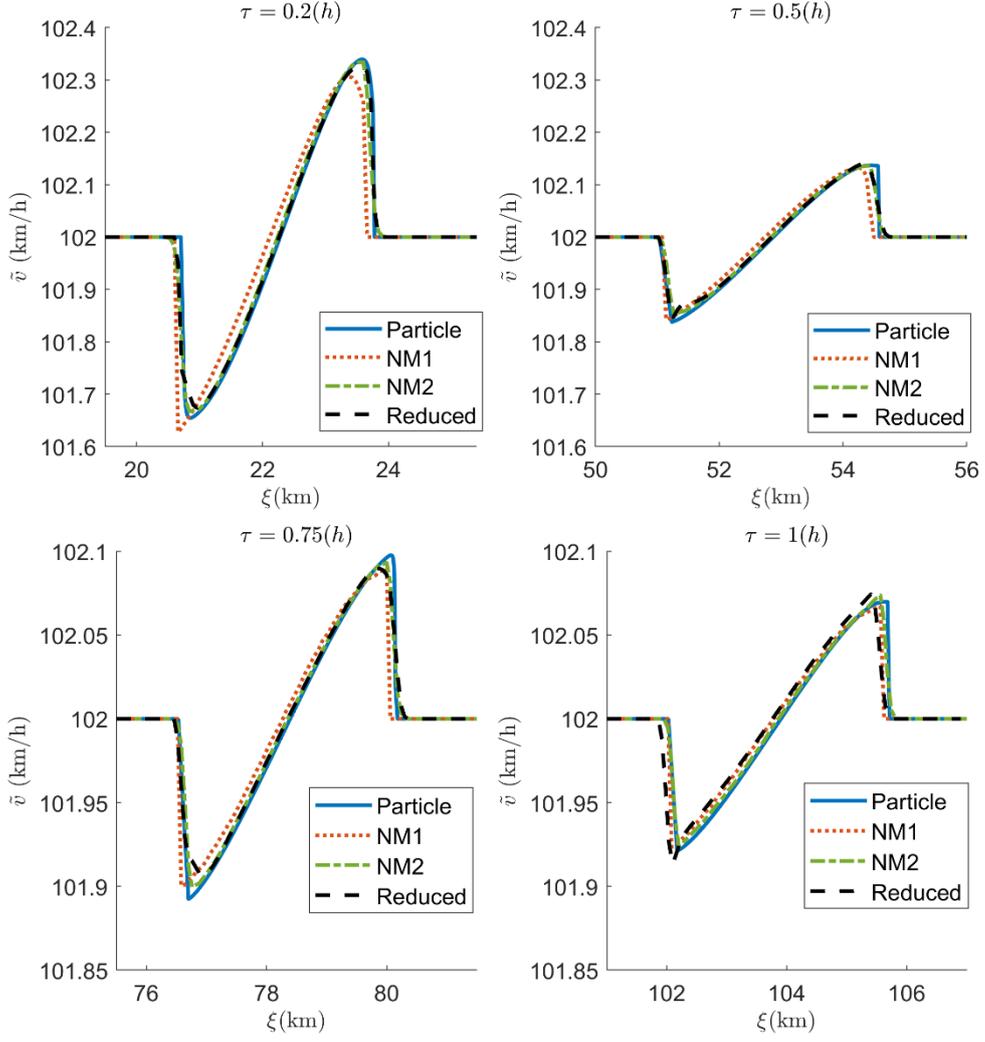

**Figure 18**: Speed profiles of the AV model (2.1)-(2.4) at $\tau = 0.2\,(h)$, $\tau = 0.5\,(h)$, $\tau = 0.75\,(h)$, and $\tau = 1\,(h)$.

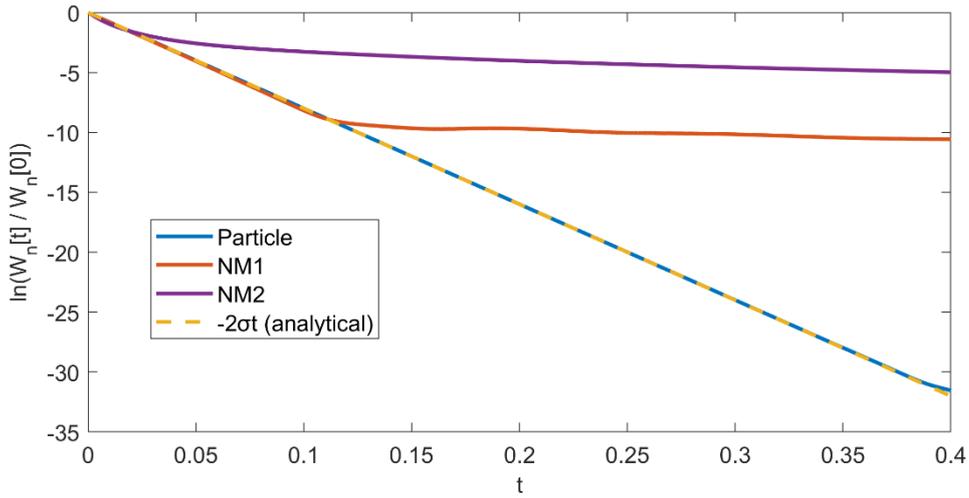

**Figure 19**: Time-evolution of $\ln(W_n^{particle}[t]/W_n[0])$, $\ln(W_n^{NM1}[t]/W_n[0])$, and $\ln(W_n^{NM2}[t]/W_n[0])$.

# 5. Conclusions

We have presented a particle method for approximating the solutions of a "fluid-like" macroscopic traffic model for automated vehicles. This method preserves all the differential inequalities that hold for the macroscopic traffic model: mass is preserved, the mechanical energy is decaying and an energy functional is also decaying. To evaluate the performance of the employed particle method, a comparison with other numerical methods for viscous compressible fluid models is provided. It is shown that although the difference between the numerical solutions is small, the proposed particle method produces sharper resolution in the numerical solutions. Moreover, the particle method under consideration has much less computational cost, in terms of CPU time, compared to the numerical methods used for viscous compressible fluid models. In addition, since the solutions of the macroscopic traffic model are approximated by the solutions of a reduced model consisting of a single nonlinear heat-type partial differential equation, the numerical solutions produced by the proposed particle method are also compared with the numerical solutions of the reduced model, verifying its high-resolution approximations. Finally, the properties of the macroscopic AV model are demonstrated by applying a traffic simulation scenario. In this scenario the macroscopic AV model is compared with the ARZ model showing that the mean flow is much higher with the use of automated vehicles and the acceleration remains close to the acceleration produced by the ARZ model.

# Acknowledgements

The research leading to these results has received funding from the European Research Council under the European Union's Horizon 2020 Research and Innovation programme/ ERC Grant Agreement n. [833915], project TrafficFluid.

# References

[1] Al-Taki, B., K. Atsou, J. Casanova, T. Goudon, P. Lafitte, F. Lagoutière, and S. Minjeaud, "Numerical Investigations of the Compressible Navier-Stokes System", *ESAIM: Proceedings and Surveys*, 70, 2021, 1-13.

[2] Ansanay-Alex, G., F. Babik, J. C. Latché, and D. Vola, "An L2-Stable Approximation of the Navier–Stokes Convection Operator for Low-Order Non-Conforming Finite Elements", *International Journal for Numerical Methods in Fluids*, 66, 2011, 555-580.

[3] Aw, A. and M. Rascle, "Resurrection of Second-Order Models of Traffic Flow", *SIAM Journal on Applied Mathematics*, 60, 2000, 916–938.

[4] Bresch, D., B. Desjardins, and D. Gérard-Varet, "On Compressible Navier–Stokes Equations with Density Dependent Viscosities in Bounded Domains", *Journal de mathématiques pures et appliquées*, 87, 2007, 227-235.

[5] Caggio, M. and D. Donatelli, "High Mach Number Limit for Korteweg Fluids with Density Dependent Viscosity", *Journal of Differential Equations*, 277, 2021, 1-37.

[6] Chertock, A. and D. Levy, "Particle Methods for Dispersive Equations", *Journal of Computational Physics*, 171, 2001, 708–730.

[7] Farjoun Y. and B. Seibold, "An Exactly Conservative Particle Method for one Dimensional Scalar Conservation Laws", *Journal of Computational Physics*, 228, 2009, 5298-5315.

[8] Gallouët, T., D. Maltese, and A. Novotny, "Error Estimates for the Implicit MAC Scheme for the Compressible Navier–Stokes Equations", *Numerische Mathematik*, 141, 2019, 495-567.

[9] Gunawan, P., "Numerical Study of Staggered Scheme for Viscous Saint-Venant Equations", *Applied Mathematical Sciences*, 8, 2014, 5349–5359.

[10] Guo, Z. and C. Zhu, "Global Weak Solutions and Asymptotic Behavior to 1D Compressible Navier–Stokes Equations with Density-Dependent Viscosity and Vacuum", *Journal of Differential Equations*, 248, 2010, 2768–2799.

[11] Haack, J., S. Jin, and J. G. Liu, "An All-Speed Asymptotic-Preserving Method for the Isentropic Euler and Navier-Stokes equations", *Communications in Computational Physics*, 12, 2012, 955-980.

[12] Hairer, E, S. P. Nørsett, and G. Wanner, *Solving Ordinary Differential Equations. I Nonstiff Problems*, Springer, Berlin, 2nd Edition, 1993.

[13] Hošek, R. and B. She, "Stability and Consistency of a Finite Difference Scheme for Compressible Viscous Isentropic Flow in Multi-Dimension", *Journal of Numerical Mathematics*, 26, 2018, 111-140.

[14] Karafyllis, I. and M. Krstic, "Global Stabilization of Compressible Flow Between Two Moving Pistons", *SIAM Journal on Control and Optimization*, 60, 2022, 1117-1142.

[15] Karafyllis, I. and M. Krstic, "Spill-Free Transfer and Stabilization of Viscous Liquid", *IEEE Transactions on Automatic Control*, 67, 2022, 4585-4597.

[16] Karafyllis, I. and M. Papageorgiou, "A Particle Method for 1-D Compressible Fluid Flow", *Studies in Applied Mathematics*, 151, 2023, 1282-1331.

[17] Karafyllis I., D. Theodosis, and M. Papageorgiou, "Constructing Artificial Traffic Fluids by Designing Cruise Controllers", *Systems & Control Letters*, 167, 2022, 105317.

[18] Karafyllis, I., Theodosis, D., and M. Papageorgiou, "Forward Completeness and Applications to Control of Automated Vehicles", submitted to *IEEE Transactions on Automatic Control* (see also arXiv:2307.11515).

[19] Karafyllis, I., F. Vokos and M. Krstic, "Feedback Stabilization of Tank-Liquid System with Robustness to Wall Friction", *ESAIM Control, Optimisation and Calculus of Variations*, 28, 2022, 81.

[20] Karafyllis, I., F. Vokos and M. Krstic, "Output-Feedback Control of Viscous Liquid-Tank System and its Numerical Approximation", *Automatica*, 149, 2023, 110827.

[21] Karper, T. K., "A Convergent FEM-DG Method for the Compressible Navier–Stokes Equations", *Numerische Mathematik*, 125, 2013, 441-510.

[22] Kazhikhov, A. V. and V. V. Shelukhin, "Unique Global Solution with Respect to Time of Initial-Boundary Value Problems for One-Dimensional Equations of a Viscous Gas", *Journal of Applied Mathematics and Mechanics*, 41, 1977, 273-282.

[23] Kotsialos A., M. Papageorgiou, C. Diakaki, Y. Pavlis and F. Middelham, "Traffic Flow Modeling of Large-scale Motorway Networks Using the Macroscopic Modeling Tool METANET," *IEEE Transactions on Intelligent Transportation Systems*, 3, 2002, 282-292.

[24] Li, T., "Global Solutions of Nonconcave Hyperbolic Conservation Laws with Relaxation Arising from Traffic Flow", *Journal of Differential Equations*, 190, 2003, 131–149.

[25] Lighthill, M. H. and G. B. Whitham, "On Kinematic Waves II: A Theory of Traffic Flow on Long Crowded Roads", *Proceedings of the Royal Society A*, 229, 1955, 317–345.

[26] Lions, P.-L., *Mathematical Topics in Fluid Dynamics, Vol.2, Compressible Models*, Oxford Science Publication, Oxford, 1998.

[27] Liotta, S. F., V. Romano, and G. Russo, "Central Schemes for Balance Laws of Relaxation Type", *SIAM Journal on Numerical Analysis*, 38, 2000, 1337-1356.

[28] Nessyahu, H., and E. Tadmor, "Non-Oscillatory Central Differencing for Hyperbolic Conservation Laws", *Journal of Computational Physics*, 87, 1990, 408-463.

[29] Nishida, T., "Equations of Motion of Compressible Viscous Fluids", *Studies in Mathematics and Its Applications*, 18, 1986, 97–128.

[30] Richards, P. I., "Shock Waves on the Highway", *Operations Research*, 4, 1956, 42–51.

[31] Shelukhin, V. V., "The Unique Solvability of the Problem of Motion of a Piston in a Viscous Gas", *Dinamika Sploshnoi Sredy*, 31, 1977, 132–150.

[32] Shelukhin V. V., "Stabilization of the Solution of a Model Problem on the Motion of a Piston in a Viscous Gas", *Dinamika Sploshnoi Sredy*, 173, 1978, 134–146.

[33] Smoller, J., *Shock Waves and Reaction-Diffusion Equations*, 2nd Edition, Springer-Verlag, New York, 1994.

[34] Theodosis, D., I. Karafyllis, G. Titakis, I. Papamichail, and M. Papageorgiou, "A Nonlinear Heat Equation Arising from Automated-Vehicle Traffic Flow Models", *Journal of Computational and Applied Mathematics*, 437, 2024, 115443.

[35] Zhang, H. M., "A Non-Equilibrium Traffic Model Devoid of Gas-like Behavior", *Transportation Research Part B*, 36, 2002, 275–290.

## Appendix: Proofs

**Proof of Proposition 1:** Using (2.43), (2.32), (2.33) and (2.24) we get:

$$\begin{aligned}
\dot{E}_n(t) &= \frac{1}{na}\sum_{i=1}^{n} H'(w_i(t))\dot{w}_i(t) + \sum_{i=2}^{n} \Phi'(nas_i(t))\dot{s}_i(t) \\
&= \frac{1}{na}\sum_{i=1}^{n} w_i(t)\beta'(w_i(t))\dot{w}_i(t) + \sum_{i=2}^{n} \Phi'(nas_i(t))(w_{i-1}(t) - w_i(t)) \\
&= \frac{1}{na} w_1(t)\beta'(w_1(t))\dot{w}_1(t) + \frac{1}{na} w_n(t)\beta'(w_n(t))\dot{w}_n(t) \\
&\quad + \frac{1}{na}\sum_{i=2}^{n-1} w_i(t)\beta'(w_i(t))\dot{w}_i(t) + \sum_{i=2}^{n} \Phi'(nas_i(t))(w_{i-1}(t) - w_i(t))
\end{aligned} \quad (A.1)$$

Moreover, using (2.34), (2.35), (2.36) we get from (A.1) that:

$$\dot{E}_n(t) = -\frac{\sigma}{na} w_1(t)\beta(w_1(t)) - w_1(t)\Phi'(nas_2(t)) + nw_1(t)K(nas_2(t))(w_2(t) - w_1(t))$$

$$-\frac{\sigma}{na} w_n(t)\beta(w_n(t)) + w_n(t)\Phi'(nas_n(t)) + nw_n(t)K(nas_n(t))(w_{n-1}(t) - w_n(t))$$

$$-\frac{\sigma}{na}\sum_{i=2}^{n-1} w_i(t)\beta(w_i(t)) + \sum_{i=2}^{n-1} w_i(t)\Phi'(nas_i(t)) - \sum_{i=2}^{n-1} w_i(t)\Phi'(nas_{i+1}(t)) \quad \text{(A.2)}$$

$$+n\sum_{i=2}^{n-1} w_i(t)K(nas_i(t))(w_{i-1}(t) - w_i(t)) + n\sum_{i=2}^{n-1} w_i(t)K(nas_{i+1}(t))(w_{i+1}(t) - w_i(t))$$

$$+\sum_{i=2}^{n} \Phi'(nas_i(t))(w_{i-1}(t) - w_i(t))$$

Since $\sum_{i=1}^{n-1} w_i(t)\Phi'(nas_{i+1}(t)) = \sum_{i=2}^{n} w_{i-1}\Phi'(nas_i(t))$, it follows by rearranging terms in (A.2) that:

$$\dot{E}_n(t) = -nw_1(t)K(nas_2(t))(w_1(t) - w_2(t))^2 - nw_n(t)K(nas_n(t))(w_{n-1}(t) - w_n(t))^2$$

$$+n\sum_{i=3}^{n-1} w_i(t)K(nas_i(t))(w_{i-1}(t) - w_i(t)) + n\sum_{i=2}^{n-2} w_i(t)K(nas_{i+1}(t))(w_{i+1}(t) - w_i(t)) \quad \text{(A.3)}$$

$$-\frac{\sigma}{na}\sum_{i=1}^{n} w_i(t)\beta(w_i(t)) = -n\sum_{i=2}^{n} K(nas_i(t))(w_{i-1}(t) - w_i(t))^2 - \frac{\sigma}{na}\sum_{i=1}^{n} w_i(t)\beta(w_i(t))$$

Recall that from (2.40), $K(s) \geq 0$ for all $s \in \left(\frac{1}{R}, +\infty\right)$. Therefore, it holds that $\sum_{i=2}^{n} K(nas_i(t))(w_{i-1}(t) - w_i(t))^2 \geq 0$. Moreover, the fact that $\sum_{i=1}^{n} w_i(t)\beta(w_i(t)) \geq 0$ is a direct consequence of the fact that $\beta(0) = 0$ and the fact that $\beta'(w) = q(w) = (1+b)^2 \frac{2b + (b-1)w}{2(b-w)^2(1+w)^2} > 0$ for all $w \in (-1, b)$, which implies that the function $\beta(w)$ defined by (2.15) is increasing. Consequently, we obtain from (A.3) that:

$$\dot{E}_n(t) = -n\sum_{i=2}^{n} K(nas_i(t))(w_{i-1}(t) - w_i(t))^2 - \frac{\sigma}{na}\sum_{i=1}^{n} w_i(t)\beta(w_i(t)) \leq 0.$$

The proof is complete. ◁

**Proof of Proposition 2:** A direct consequence of Proposition 1 is the following inequality:

$$E_n(t) = \frac{1}{na}\sum_{i=1}^{n} H(w_i(t)) + \frac{1}{na}\sum_{i=2}^{n} \Phi(nas_i(t))$$

$$\leq E_n(0) = \frac{1}{na}\sum_{i=1}^{n} H(w_i(0)) + \frac{1}{na}\sum_{i=2}^{n} \Phi(nas_i(0)) \quad \text{(A.4)}$$

as long as the solution of (2.32)-(2.36) is defined. Using the fact that $H(w) \geq 0$ for all $w \in (-1, b)$ (see (2.24)) and the fact that $\Phi(s) \geq 0$ for all $s \in \left(\frac{1}{R}, +\infty\right)$ (a consequence of (2.40)), we get from (A.4) that:

$$\Phi(nas_i(t)) \leq \sum_{j=1}^{n} H(w_j(0)) + \sum_{j=2}^{n} \Phi(nas_j(0)), \text{ for } i = 2,...,n \quad (A.5)$$

$$H(w_i(t)) \leq \sum_{j=1}^{n} H(w_j(0)) + \sum_{j=2}^{n} \Phi(nas_j(0)), \text{ for } i = 1,...,n \quad (A.6)$$

Using the fact that $\lim_{w \to (-1)^+} H(w) = \lim_{w \to b^-} H(w) = +\infty$ and (A.6) we conclude that there exist constants $-1 < w_{\min} < w_{\max} < b$ such that $w_i(t) \in [w_{\min}, w_{\max}]$ for $i = 1,...,n$ as long as the solution of (2.32)-(2.36) is defined. Using (2.41) and (A.5) we also conclude that there exists a constant $s_{\min} > \frac{1}{naR}$ such that $s_i(t) \geq s_{\min}$ for $i = 2,...,n$ as long as the solution of (2.32)-(2.36) exists. Notice that boundedness of $w_i(t) \in [w_{\min}, w_{\max}]$, for $i = 1,...,n$ as long as the solution of (2.32)-(2.36) is defined, implies that $\dot{s}_i(t)$ are bounded for all $i = 2,...,n$. Moreover, by virtue of (2.40) and since $s_i(t) \geq s_{\min}$ we also obtain boundedness of $\Phi'(s_i)$ and $K(s_i)$, which in conjunction with continuity of $\beta'(w)$ and the fact that $w_i(t) \in [w_{\min}, w_{\max}]$, for $i = 1,...,n$, as long as the solution of (2.32)-(2.36) is defined, implies that $\dot{w}_i(t)$, are bounded for all $i = 1,...,n$. The latter, in conjunction with (2.46) of Proposition 1, and Theorem 1 in [18], implies that every solution of (2.32)-(2.36) is defined for all $t \geq 0$. ◁

**Proof of Proposition 3:** Using (2.42) we obtain

$$\Phi''(s) = \frac{\sigma}{s^2} \kappa\left(\frac{1}{s}\right)$$

which due to definition (2.38) implies that

$$\Phi''(s) = \frac{\sigma}{a} K(s) \quad (A.7)$$

Combining definition (2.45) and (A.7) we get that for $i = 2,...,n-1$

$$\dot{\varphi}_i(t) = \beta'(w_i(t))\dot{w}_i(t) + n^2 aK(nas_{i+1}(t))(w_i(t) - w_{i+1}(t)) \\ - n^2 aK(nas_i(t))(w_{i-1}(t) - w_i(t)) \quad (A.8)$$

Next, using (2.35), (2.44), (2.45), and (A.8) we get:

$$\dot{W}_n(t) = \frac{1}{na} \sum_{i=2}^{n-1} \varphi_i(t) \dot{\varphi}_i(t)$$

$$= \frac{1}{na} \sum_{i=2}^{n-1} \varphi_i(t) \beta'(w_i(t)) \dot{w}_i(t) + n \sum_{i=2}^{n-1} \varphi_i(t) K(nas_{i+1}(t))(w_i(t) - w_{i+1}(t))$$

$$- n \sum_{i=2}^{n-1} \varphi_i(t) K(nas_i(t))(w_{i-1}(t) - w_i(t))$$

$$= \frac{1}{na} \sum_{i=2}^{n-1} \varphi_i(t) \Big( -\sigma \beta(w_i(t)) + na\, \Phi'(nas_i(t)) - na\, \Phi'(nas_{i+1}(t)) \quad\quad\quad (A.9)$$

$$+ n^2 a K(nas_i(t))(w_{i-1}(t) - w_i(t)) + n^2 a K(nas_{i+1}(t))(w_{i+1}(t) - w_i(t)) \Big)$$

$$+ n \sum_{i=2}^{n-1} \varphi_i(t) K(nas_{i+1}(t))(w_i(t) - w_{i+1}(t)) - n \sum_{i=2}^{n-1} \varphi_i(t) K(nas_i(t))(w_{i-1}(t) - w_i(t))$$

$$= -\frac{\sigma}{na} \sum_{i=2}^{n-1} \varphi_i(t) \Big( \beta(w_i(t)) + \sigma^{-1} na \big( \Phi'(nas_{i+1}(t)) - \Phi'(nas_i(t)) \big) \Big)$$

Using again definition (2.45), we get from (A.9) that

$$\dot{W}_n(t) = -\frac{\sigma}{na} \sum_{i=2}^{n-1} \varphi_i^2(t) \leq 0.$$

This completes the proof. ◁